\documentclass[12pt,hidelinks]{article}

\usepackage[utf8]{inputenc}
\usepackage[T1]{fontenc}

\usepackage[colorlinks=false]{hyperref}

\usepackage{amsmath, amsfonts,amsthm,amssymb}
\usepackage{bbm}

\usepackage{graphicx}

\usepackage{txfonts}

\usepackage{setspace}

\usepackage{pgf,tikz}
\usepackage{pgfplots}
\usetikzlibrary{spy}
\usetikzlibrary{backgrounds}
\usetikzlibrary{decorations}

\usepackage[font=footnotesize,labelfont=bf]{caption}

\usepackage{color}

\usepackage{customcommands}

\theoremstyle{theorem}

\newtheorem{theorem}{Theorem}
\newtheorem{lemma}[theorem]{Lemma}

\theoremstyle{definition}
\theoremstyle{remark}
\newtheorem{remark}{Remark}

\usepackage[sort&compress]{natbib}
\usepackage{snapshot}

\usepackage{booktabs}
\usepackage{multicol}

\usepackage{longtable}

\usepackage[boxed, ruled, lined]{algorithm2e}

\renewcommand{\argmin}{\operatorname*{argmin}}

\usepackage{subcaption}

\usepackage{enumerate}
\usepackage{csquotes}

\usepackage{changes}
\definechangesauthor[name=MS, color=orange]{MS}

\title{
Smoothing splines for discontinuous signals 
}
\author{Martin Storath\footnote{Lab for Mathematical Methods in Computer Vision and Machine Learning, Technische Hochschule Würzburg-Schweinfurt, Schweinfurt, Germany.}, Andreas Weinmann\footnote{Department of Mathematics and Natural Sciences, Hochschule Darmstadt, Darmstadt, Germany. 
} }

\begin{document}

\maketitle

\begin{abstract}
Smoothing splines are twice differentiable by construction, so they cannot capture potential discontinuities in the underlying signal. 
In this work, we consider a special case
 of the weak rod model of Blake and Zisserman (1987) that allows for discontinuities penalizing their number by a linear term.
    The corresponding estimates are cubic smoothing splines with discontinuities (CSSD) which serve as representations of piecewise smooth signals
 and facilitate exploratory data analysis. 
However, computing the estimates requires solving a non-convex optimization problem.
So far, efficient and exact solvers  exist only for a discrete approximation based on  equidistantly sampled data.
In this work, we propose an efficient solver for the continuous minimization problem with non-equidistantly sampled data. Its worst case complexity is quadratic in the number of data points, and if the number of detected discontinuities scales linearly with the signal length, we observe linear growth in runtime. 
This efficient algorithm allows to use cross validation for automatic  selection of the hyperparameters within a reasonable time frame on standard hardware. We provide a reference implementation and supplementary material. We demonstrate the applicability of the approach for  the aforementioned tasks   using both simulated and real data. 
\end{abstract}

\maketitle
\section{Introduction}

Assume that we are given the approximate values $y_i = g(x_i) + \epsilon_i$ 
of a function $g$ at data sites $x_1, \ldots,x_N,$
and an estimate $\delta_i$  of the standard deviation of the errors $\epsilon_i$
which are uncorrelated with zero mean.
If we can assume that $g$ is a smooth function,
then we may try to estimate $g$ using a (cubic) smoothing spline,
a flexible and widely-applicable approach to curve estimation \citep{silverman1985some, hastie2009elements}.
Using the notations of \cite{de1978practical}, a cubic smoothing spline is the solution $\hat f$ of the variational problem
\begin{equation}\label{eq:spline}
\min_{f \in C^2}	p \sum_{i=1}^N \left(\frac{y_i - f(x_i)}{\delta_i}\right)^2   +  (1-p) \int_{x_1}^ {x_N}   (f''(t))^2 \,dt.
\end{equation}
Here, the minimum is taken over all twice continuously differentiable functions on $[x_1, x_N]$. 
The minimizer is a compromise between the conflicting goals smoothness, measured by the squared euclidean norm of the derivative, and closeness to data. The \enquote{stiffness} parameter $p \in (0,1)$ controls the relative weight of the two goals.
It is well known that the minimizer  $\hat f$  of \eqref{eq:spline} is a {cubic spline,} that is,  $\hat f$ is a piecewise cubic polynomial function which is twice continuously differentiable \citep{de1978practical}. 
If the underlying function $g$ has {discontinuities} (or {breaks/jumps}), or in other words, $g$ is only {piecewise smooth}, the classical smoothing spline cannot capture these discontinuities.

Signals with discontinuities 
appear in numerous applications; for example, the cross-hy\-bridization of DNA \citep{snijders2001assembly,drobyshev2003specificity,hupe2004analysis}, 
the reconstruction of brain sti\-muli \citep{winkler2005don}, 
single-molecule analysis \citep{joo2008advances,loeff2021autostepfinder},
cellular ion channel functionalities \citep{hotz2013idealizing},
photo-emission spectroscopy \citep{frick2014multiscale}
and the rotations 
of the bacterial flagellar motor 
\citep{sowa2005direct}; see also 
\cite{kleinberg2006algorithm} and \cite{frick2014multiscale} for further examples.

\begin{figure}[!t]
	\centering
	\includegraphics[width=0.9\textwidth]{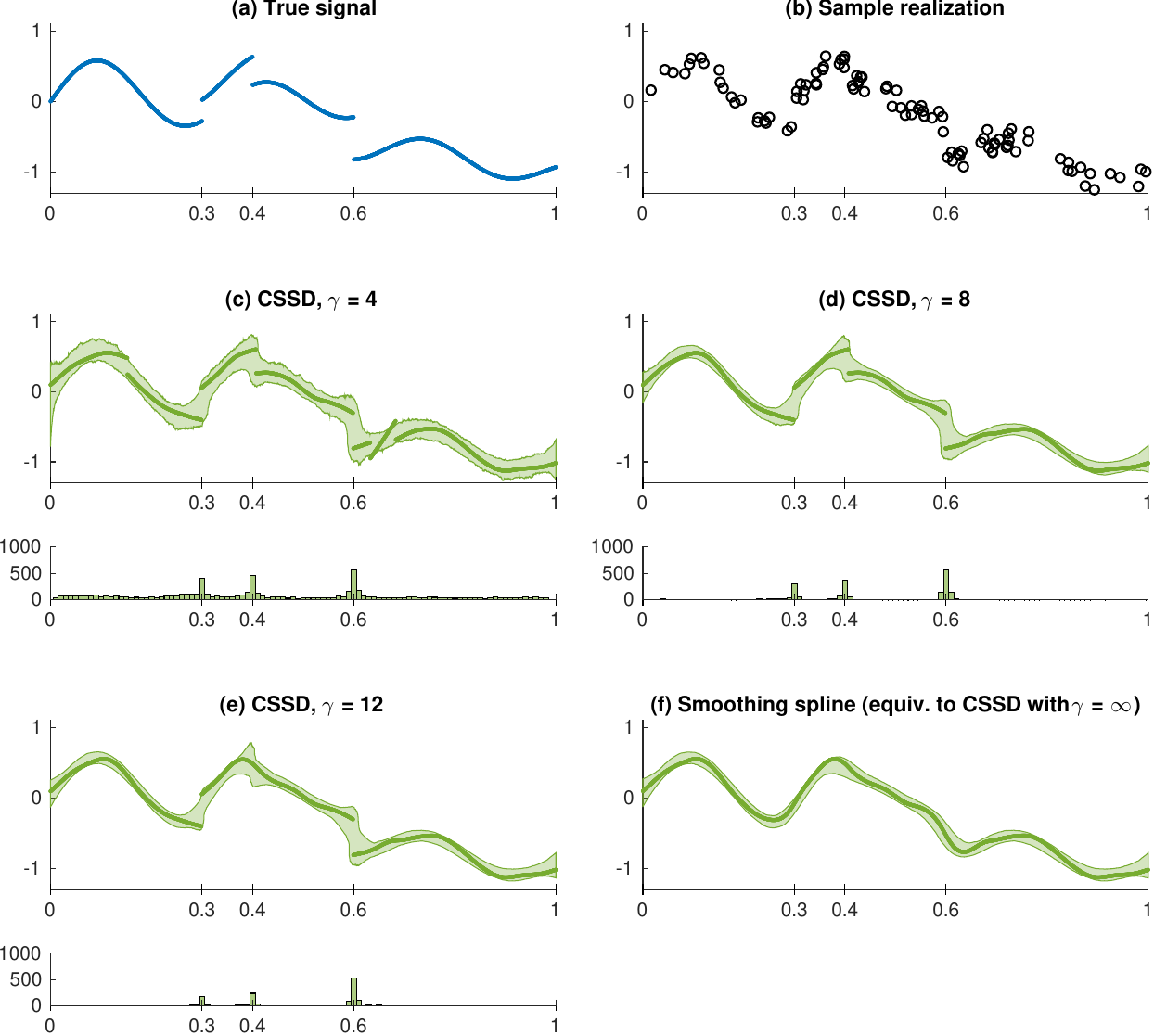}
	\caption{A synthetic signal is sampled at $N = 100$ random data sites $x_i$ 
	and corrupted by zero mean Gaussian noise with standard deviation 
	$0.1.$
	The results of the discussed model are shown for $p=0.999$ and different parameters of $\gamma,$ where $\gamma=\infty$  corresponds to classical smoothing splines.
	The thick lines 
represent the results of the shown sample realization.	The x-ticks indicate the true discontinuity locations. The shaded areas depict the $2.5 \%$ to $97.5 \%$  (pointwise) quantiles of $1000$ realizations. The histograms under the plots show the frequency of the detected discontinuity locations over all realizations.}
	\label{fig:synthetic}
\end{figure}

We are interested in the case where the locations of the discontinuities are {unknown}. 
(If the locations of the discontinuities are known a priori, one simply may compute smoothing splines on the intervals between two points of discontinuities.)
In their landmark work, \cite{blake1987visual} proposed a variational model for this task based on piecewise regression with smoothing splines and linear penalties 
on the number of jumps and creases, called the \textit{weak rod model.}
Here, we study the special case of the weak rod model without creases, and for
 discrete, potentially non-equidistant, sampling points. Using the above notation of smoothing splines, it can be formulated as:
\begin{equation}\label{eq:pcw_spline_lagrange}
\min_{f, J} p \sum_{i=1}^N \left(\frac{y_i - f(x_i)}{\delta_i}\right)^2   +  (1-p) \int_{[x_1, x_N] \setminus J}   (f''(t))^2 \,dt	 + \gamma |J|. 
\end{equation}
Here the minimum is taken over all possible sets of discontinuities between two data sites $J \subset [x_1, x_N]\setminus \{x_1, \ldots, x_N\}$
 and all functions $f$ that are twice continuously differentiable away from the discontinuities.
The last term is a penalty for the number of discontinuities $|J|$
weighted by a parameter $\gamma > 0.$ 
In \eqref{eq:pcw_spline_lagrange}, we search for 
a global minimizer which consists of a discontinuity set $\hat{J}$ and a piecewise cubic spline $\hat{f}$ with discontinuities at $\hat{J}.$
The solution $\hat{f}$ of the specific minimization problem \eqref{eq:pcw_spline_lagrange} is a \textit{cubic smoothing spline with discontinuities} which we abbreviate as \textit{CSSD.}
A first example of a CSSD is provided in Figure~\ref{fig:synthetic}.
A CSSD may be used for exploratory data analysis, 
as predictor, and as function estimator; just like a classical smoothing spline  \citep{silverman1985some}, but with the  extension that the underlying function may have discontinuities. Additionally, a detected discontinuity can be seen as a type of changepoint.

The considered instance of the weak rod model \eqref{eq:pcw_spline_lagrange} is interesting and useful on its own, because it comprises two widely used regression models: For sufficiently large~$\gamma,$ a CSSD coincides with the $C^2$ continuous cubic smoothing spline. In the limit $p \to 0,$ a CSSD tends to a piecewise linear regression function,
which were studied e.g. in \cite{kleinberg2006algorithm}, \cite{friedrich2008complexity},  \cite{storath2019smoothing}.
If the second order derivative is replaced by the first order derivative in
\eqref{eq:pcw_spline_lagrange}, we obtain the Mumford-Shah model \citep{mumford1989optimal}, and the according limit $p \to 0$
leads to a piecewise constant regression model, also referred to as Potts model \citep{winkler2002smoothers, winkler2003image},
and studied in a series of further works, including the works of \cite{jackson2005algorithm},  \cite{boysen2009consistencies}, and \cite{killick2012optimal}, to mention only a few.

The literature describes  algorithms
for a variant of \eqref{eq:pcw_spline_lagrange} with a discretized roughness penalty for equidistantly sampled data sites, i.e. $x_{i+1} - x_{i} = const$ for all $i.$ 
 For solving this discrete problem, \cite{blake1987visual} proposed a graduated non-convexity algorithm, an algorithm based on Hopfield’s neural model and a Viterbi-type algorithm but these algorithms are not guaranteed to yield a global minimum in general. Straightforward adaption of the slightly different dynamic programming approaches of \cite{blake1989comparison} or \cite{winkler2002smoothers}, one may obtain a global minimizer in cubic worst case time complexity.  \cite{storath2019smoothing} proposed an algorithm which solves the discretized problem for equidistantly sampled data in quadratic complexity.

However, restriction to equidistantly sampled data means a major limitation
because \textit{(i)} many types of data can be acquired only with varying distances of the measurements sites,
and \textit{(ii)} even when the data sites are equidistantly sampled,
non-equidistant data points occur naturally when performing $K$-fold cross validation for parameter selection.
The present formulation with continuous roughness penalty \eqref{eq:pcw_spline_lagrange} naturally deals  with non-equidistant data sites $x_i,$ but solving it is more challenging than the  discrete equidistant setting.
 There is -- to the authors knowledge -- no dedicated paper presenting an exact solver for \eqref{eq:pcw_spline_lagrange}. 
After the reformulation to a reduced search space (see Section \ref{sec:reduction}), one may combine methods of dynamic programming  \citep{blake1989comparison,winkler2002smoothers, killick2012optimal}  and smoothing splines  \citep{reinsch1967smoothing},  to obtain an exact solver for \eqref{eq:pcw_spline_lagrange}.
However, such a \enquote{baseline} solver has worst case cubic time complexity, and hence applying it is costly even for problems of moderate size. This hampers in particular usage of standard parameter selection strategies, such as cross validation, since they require solving \eqref{eq:pcw_spline_lagrange} for numerous model parameters $\gamma$ and $p.$ 
The present paper develops an efficient exact solver for \eqref{eq:pcw_spline_lagrange}.
Details on the contribution of this paper follow after a further discussion of  prior and related work.

\subsection{Prior and related work}
\cite{blake1987visual} have introduced the weak rod model in the context of image processing
along with a two-dimensional generalization which they called the weak plate model. 
Equation \eqref{eq:pcw_spline_lagrange} is a special case of the weak rod model 
in the sense that it does not allow creases (corresponding to the weak rod model with a sufficiently high crease penalty), and
that it is formulated with a potentially weighted discrete data term which naturally appears  when working with sampled data.
\cite{blake1987visual} utilized the
penalty factors $\mu, \alpha > 0$ which 
are related to the factors $p, \gamma$  in \eqref{eq:pcw_spline_lagrange} by  
$\mu = \sqrt[4]{(1-p)/p}$ and $\alpha = \gamma/p.$
(The present paper adopts the $p$-parametrization used in a standard reference on splines by  \cite{de1978practical}.)
\cite{blake1987visual} have obtained results on scale and sensitivity in discontinuity detection for the fully continuous version.
In particular, they derived a contrast threshold  referring to the minimum detectable jump height for an isolated bi-infinite step function in dependance of the model parameters which is given by $h_0 = 2^{\frac{3}{4}}\sqrt{\alpha / \mu}.$

If we replace the second order derivative with the first order derivative in Equation \eqref{eq:pcw_spline_lagrange},
we obtain the weak string model which coincides with the one-dimensional case of the famous Mumford-Shah model 
for edge preserving smoothing \citep{mumford1985boundary,mumford1989optimal}.
This first order model has been investigated in more detail in the literature than the second order model; in particular, an exact  algorithm of cubic complexity for the discretized first order model has been proposed in an early work by \cite{blake1989comparison}.
\cite{blake1987visual} have described a fundamental limitation of the first order model termed the gradient limit; it refers to the undesired introduction of spurious discontinuities when the gradient exceeds a threshold.
The second order model considered here does not have this limitation.

The present model is closely connected to smoothing splines which were developed in the works of \cite{schoenberg1964spline} and \cite{reinsch1967smoothing}, and the basic idea can be tracked back to the work of \cite{whittaker1923new}. \cite{silverman1984spline} has shown that spline smoothing
approximately corresponds  to smoothing by a kernel method with bandwidth depending on the local density of the data sites.
Comprehensive treatments on smoothing splines are given in the books by \cite{de1978practical}, by  \cite{wahba1990spline}, and by \cite{green1993nonparametric}, as well as in the paper by \cite{silverman1985some}. 
(In the present paper, we use the notations and conventions of \cite{de1978practical}.)
The signal and image processing point of view is discussed in the papers of \cite{unser1999splines, unser2002splines}.
Smoothing splines have become a standard method for statistical data processing, and they are discussed in standard literature on the topic \citep{hastie2009elements}.
Recent contributions on splines in a statistical context are generalizations for Riemannian manifolds \citep{kim2021smoothing}, locally adaptive splines for Bayesian model selection in additive partial linear models \citep{jeong2021bayesian}, knot estimation for linear splines, \citep{yang2021estimation}, and efficiency improvements  for splines with multiple predictors \citep{meng2022smoothing}, to mention only a few.

Besides the above mentioned models of spline-type, 
many other piecewise regression models with discontinuity or jump penalty
have been used. They typically lead to combinatorial optimization problems,
and fast solvers are a central question \citep{winkler2002smoothers, killick2012optimal,frick2014multiscale}.
Efficient algorithms for piecewise constant or piecewise polynomial regression functions have been proposed by \cite{auger1989algorithms}, \cite{winkler2002smoothers}, \cite{jackson2005algorithm}, \cite{friedrich2008complexity},  \cite{little2011generalized2}, and \cite{killick2012optimal}, and extended to indirectly measured data \citep{weinmann2015iterative,storath2014jump} and to non-linear  data spaces \citep{storath2017jumppenalized, weinmann2016mumford}. Parallelized versions have been considered by \cite{tickle2020parallelization}.
For the piecewise constant case, \cite{boysen2009consistencies} have obtained consistency and convergence rates of the estimates in $L^2([0,1)).$ 
A detailed treatment of piecewise constant regression is given by \cite{little2011generalized,little2011generalized2}.

Regarding automatic selection of the model parameters, there are some common approaches for related piecewise regression models;
for example information based criteria \citep{zhang2007modified, yao1988estimating},  an interval criterion \citep{winkler2005don}
and several variants of cross validation \citep{arlot2011segmentation}.
For classical smoothing splines,
generalized cross validation \citep{silverman1985some,craven1979smoothing,golub1979generalized} is frequently used. The authors are not aware of an automatic parameter selection strategy specifically developed for the considered model \eqref{eq:pcw_spline_lagrange}.

There are also conceptually different  methods for estimating discontinuous regression functions of spline type. 
One such approach was proposed by \cite{koo1997spline} where knots are stepwise refined by knot addition, basis deletion, knot-merging, and the final model is determined by the Bayes information criterion. 
Another possibility is a changepoint based approach:
one might first use an existing method for changepoint detection to determine discontinuity locations; for example the CUSUM method \citep{page1954continuous}, Bayesian changepoint inference \citep{fearnhead2006exact}, wild binary segmentation \citep{fryzlewicz2014wild}, a narrowest over the threshold method \citep{baranowski2019narrowest}, a
Bayesian ensemble approach  \citep{zhao2019detecting}, nonparametric maximum likelihood approaches \citep{zou2014nonparametric, haynes2017computationally}, or multiscale testing \citep{frick2014multiscale} could be used.
 (See \cite{paper:an-evaluation-of-change-point-detection-algorithms} and the references therein for an overview and a comparison of selected changepoint detection methods.)
Then, splines can be fitted between two detected changepoints; either smoothing splines directly to the data or interpolating splines to the corresponding signal estimates. We point out that such a two stage approach emphasizes the changepoint detection aspect, and as splines are typically not involved in the detection process, they are not necessarily the \enquote{natural} models for the data between two changepoints.

For denoising non-smooth signals, shrinkage of wavelet coefficients is frequently used \citep{donoho1994ideal}. 
The book by \cite{mallat2008wavelet} provides an overview on wavelet shrinkage.
By their multiscale subdivision algorithms, wavelet methods  are computationally extremely efficient.
In contrast to the model considered in this work, they typically rely on equidistant sampling points and do not result in piecewise smooth regression functions. 

\subsection{Contribution}

We first discuss  basic properties of a CSSD,
that is, the properties of the minimizers of \eqref{eq:pcw_spline_lagrange}.
In particular, in view of the well-definedness of the estimator, we provide  uniqueness results for the optimization problem with respect to both function evaluations and partitions in an almost everywhere sense.

The main contribution of this paper is an efficient algorithm for computing 
a solution of the problem \eqref{eq:pcw_spline_lagrange}, meaning a global minimum of the target function.
The algorithm is developed in three steps:
\textit{(i)} We show that we may restrict the search space for the discontinuity set $J$ to  the midpoints of the data sites, and we use dynamic programming to reduce the number of possible configurations.
\textit{(ii)} We propose a procedure that computes the spline energies
for a signal of increasing length in constant time per new element. It allows to obtain the required spline energies on all (discrete) intervals efficiently.
We point out that using the smoothing spline algorithms of \cite{reinsch1967smoothing} and of \cite{de1986efficient} lead to significantly higher computational costs for this specific task.
\textit{(iii)} We show that the computation of the spline energies is compatible with  two different pruning strategies, the PELT pruning of \cite{killick2012optimal} and the FPVI pruning of \cite{storath2014fast}.
The worst case time complexity of the proposed algorithm is $O(N^2)$ where $N$ is the number of data points. If the number of detected discontinuities scales linearly with $N,$ we observe linear scaling of the runtime. 

We provide a ready-to-use reference implementation in Matlab, available for download on Github.\footnote{\url{https://github.com/mstorath/CSSD}}.
The novel solver is notably faster than a baseline solver which relies on standard Python modules and which does not use of the methods developed in the present work.
We further implement a strategy for  selecting the model parameters ($p$ and $\gamma$) automatically  based on $K$-fold cross validation.
The proposed fast algorithm allows to perform this selection strategy within a reasonable time frame on standard hardware.
Numerical examples with synthetic and real data demonstrate the potential of CSSD
as function estimator for discontinuous signals, as basis for a changepoint detector,
 and as a tool for exploratory data analysis.

\section{Efficient computation of a CSSD and uniqueness result}

\subsection{Reformulation and basic properties of the solutions}\label{sec:basic}

Throughout this paper, we assume that the data sites satisfy
$x_1 < x_2 < \ldots < x_N.$ 
If the data does not satisfy this constraint, 
we may merge data sites into a single data point by weighted averaging over the $y$-values
of coinciding $x$-values; see e.g. \cite{hutchinson1986algorithm}.

Let $J \subset [x_1, x_N]\setminus \{x_1, \ldots, x_N\}$ be a set of discontinuities 
of size $|J| = K.$ It is convenient to sort the elements of $J$ ascendingly
so that we can access the elements by an index, e.g. $J_1$ for the smallest element in $J.$
Its complement in $[x_1, x_N]$ consists of $K+1$ half open and open intervals, denoted by  $I_1, \ldots I_{K+1}$
\[
	[x_1, x_N] \setminus J = I_1 \cup \ldots \cup I_{K+1}.
\]
We denote the (ordered) set of these intervals by $\Pc(J)$
\[
	\Pc(J) = \{I_1, \ldots, I_{K+1} \}.
\]
The minimization problem \eqref{eq:pcw_spline_lagrange} now can be reformulated equivalently in terms of the discontinuity set and the corresponding intervals as
\begin{equation}\label{eq:pcw_spline}
\min_{J \subset [x_1, x_N]\setminus \{x_1, \ldots, x_N\}} \sum_{I \in \Pc(J)} \Ec_I	 + \gamma |J|
\end{equation}
where $\Ec_I$ is given by
\begin{equation}\label{eq:epsilon_I}
	\Ec_I =  \min_{f \in C^2(I)} p \sum_{i: x_i \in I} \left(\frac{y_i - f(x_i)}{\delta_i}\right)^2   +  (1-p) \int_{I}   (f''(t))^2 \,dt.
\end{equation}
As in \eqref{eq:spline}, the model parameter $p \in (0, 1)$ controls the relative weight of the smoothness term and the data fidelity term, and $\delta_i$ is an estimate of the standard deviation of the errors at data site $x_i.$
Note that \eqref{eq:epsilon_I} describes the minimum energy of a classical smoothing spline on the interval $I.$

A standard procedure to obtain $\Ec_I$ consists in computing a minimizer $f_I$ of the functional in \eqref{eq:epsilon_I} using the algorithm of \cite{reinsch1967smoothing} on the data in the interval $I,$ and evaluating \eqref{eq:epsilon_I} with that $f_I$.  
The computational costs of this approach are linear in the number of data points falling in the interval $I.$  
A key element of the fast algorithm we are going to propose in the next section is that it bypasses the costly explicit computation of the $f_I,$ and computes $\Ec_I$ directly using a recurrence relation.

From \eqref{eq:pcw_spline} we see that every discontinuity set $J$
has a unique piecewise cubic spline function $f_J$ which is defined for all $I \in \Pc(J)$ by the unique minimizer of \eqref{eq:epsilon_I}. So once an optimal discontinuity set is known, 
the corresponding $f_J$ can be computed by a standard method such as Reinsch's algorithm. 
In opposite direction, $J$ consists of the discontinuities of $f$ and the discontinuities of $f'.$ 

Next we state two basic properties of an optimal discontinuity set $\hat{J}.$ 
\begin{lemma}\label{lem:J_star}
Let $\hat{J}$ be a solution of the minimization problem \eqref{eq:pcw_spline}. 
\begin{enumerate}
\item $\hat{J}$ has at most one element between two adjacent data sites and
\item $\hat{J}$  has at most $\ceil{N/2} - 1$ elements.  
\end{enumerate}

\end{lemma}
The proof is given in the supplementary material.

We develop an efficient solver for \eqref{eq:pcw_spline} -- the equivalent formulation of \eqref{eq:pcw_spline_lagrange} -- in three steps:
\textit{(i)} reduction of the search space for the discontinuities to a discrete set and reducing  the number of configurations using dynamic programming,
\textit{(ii)}  efficient update scheme for the necessary spline energies,
\textit{(iii)} search space pruning which is compatible with the order of computation of the energy updates.

\subsection{Reduction of the search space for the discontinuity set}\label{sec:reduction}

\paragraph*{Reducing the discontinuity set to midpoints of data sites.}
It is a basic property of the classical smoothing spline that its energy $\Ec_I$ (defined in \eqref{eq:epsilon_I})
only depends on the data sites that are contained in the interval $I.$ The reason for this is that the spline can be extended linearly beyond the extremal data sites without increasing the energy; see e.g. \cite{silverman1985some}.
So if $I$ and $I'$ are two intervals
which contain the same data sites, that is, if  $\{x_i: x_i \in I\} = \{x_i: x_i \in I'\},$ then 
	$
	\Ec_I = \Ec_{I'}.
	$
As a direct consequence, a shift of a discontinuity location 
between two data sites $x_i, x_{i+1}$ does not change the functional value \eqref{eq:pcw_spline}.
So a discontinuity  may be located at an arbitrary position between two data points.
Having noticed this, 
we consider two discontinuity sets as equivalent,
if all their respective elements lie between the same data sites.
\begin{figure}
	\includegraphics[width=1\textwidth]{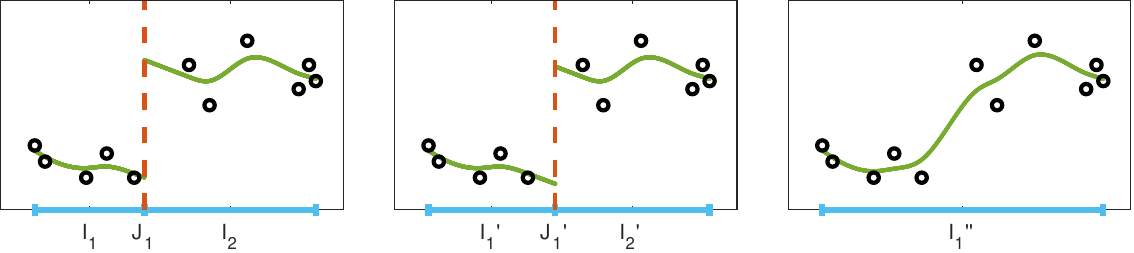}
	\caption{A toy data set (circles) and three possible discontinuity configurations (visualized on abscissa) along with their corresponding piecewise splines (solid curves) are shown ($x$ on abscissa, $y$ on ordinate). 
	The first two configurations give the same function value $\Ec_{I_1} + \Ec_{I_2} + \gamma = \Ec_{I_1'} + \Ec_{I_2'} + \gamma,$ when plugged into the target function in \eqref{eq:pcw_spline_lagrange}. Hence, these  configurations are equally good in the sense of the model \eqref{eq:pcw_spline_lagrange}. So we consider the discontinuity locations $J_1$ and $J_1'$ as equivalent and use the midpoint (here $J_1',$ middle tile) as natural representative. The third tile shows a configuration with no discontinuity which yields the	target function value $\Ec_{I_1''}.$ 
If $\gamma$ is sufficiently small ($\gamma < \Ec_{I_1''} - \Ec_{I_1'} - \Ec_{I_2'}$) the   configurations with discontinuity have a better target function value than the one without discontinuity.}
	\label{fig:jumpLocations}
\end{figure}
The midpoint between two data sites 
is a natural choice for a representative of the equivalence class.
This is illustrated in Figure~\ref{fig:jumpLocations} using a toy example.
Hence,  it is sufficient to consider the discrete set of midpoints $\Mc$  of the data sites, 
\begin{equation}\label{eq:midpoints}
\Mc := \{\tfrac{1}{2}{(x_i+x_{i+1})}: i = 1, \cdots, N-1\},
\end{equation}
as search space for $J.$  
Thus, we have reduced the optimization problem \eqref{eq:pcw_spline} to
the discrete optimization problem
\begin{equation}\label{eq:pcw_spline_discrete}
\min_{J \subset \Mc} \sum_{I \in \Pc(J)} \Ec_I	 + \gamma |J|.
\end{equation}
Because $\Mc$ is a finite set it follows that the minimization problem \eqref{eq:pcw_spline_discrete}, and so the problems \eqref{eq:pcw_spline_lagrange} and \eqref{eq:pcw_spline}, indeed have  a minimizer.

\begin{remark}
For the algorithm, which we are developing here, the concrete representative of the  discontinuity set is irrelevant.
The choice of the representative comes into play when displaying the final result. 
For this purpose, it  seems natural to take the midpoints.
Yet, there might be scenarios where other points -- or even the entire interval -- are reasonable representatives for the discontinuity locations.
(One may think of a similar situation when dealing with the ordinary median which can be defined as 
 $\argmin_\mu |y_i - \mu|.$
 For data of even length, all values between the two central data points satisfy the minimality property, and most commonly, the mean value between the central points is chosen. However, for implementing the method of \cite{weinmann2014l1potts}, considering the upper and lower points as median turned out to be useful.)
\end{remark}

\paragraph*{Reducing the number of configurations by dynamic programming.}
The reduced form \eqref{eq:pcw_spline_discrete} constitutes a discrete one-di\-men\-sional partitioning problem which can be solved by dynamic programming.
The approach works similar to the corresponding algorithms on related partitioning problems;
see  \cite{friedrich2008complexity, kleinberg2006algorithm, winkler2002smoothers,jackson2005algorithm}.
For completeness we briefly describe the procedure.
For the formulation it is convenient, to identify an interval $I$ with the indices of the contained data sites.
For example, the interval that contains the data sites $x_l, x_{l+1}, \ldots, x_{r}$
is identified with the (ordered) set of indices $\{l, l+1, \ldots, r\},$ 
or abbreviated in Matlab-type notation $\{ l:r\}.$ 
In particular, we write
$	\Ec_{\{ l:r\}}$ for $\Ec_I.$

We now consider the functional in \eqref{eq:pcw_spline_discrete} 
 on a reduced data set $(x_1,y_1),$ $\ldots,$ $(x_r, y_r)$ for $r \leq n$
 and denote it by $F_r;$ that is 
$
	F_r(J) =  \sum_{I \in \Pc(J)} \Ec_{I} + \gamma |J|.
$
where $J \subset \Mc \cap [x_1,x_r]$ and $\Ec_{I}$ is only evaluated for the data points $x_1, \ldots, x_r.$
The minimal functional value of $F_r,$ denoted by $F^*_r = \min_{J}  F_r(J)$,
satisfies the Bellman equation
\begin{equation}\label{eq:recurrencePenalized}
    F^*_r 
    =  \min_{l= 1, \ldots, r} \Big\{  \Ec_{\{ l:r\}} + \gamma + F^*_{l-1} \Big\}
\end{equation}
where we let $F^*_0 = -\gamma,$ see \cite{friedrich2008complexity}.
As $\Ec_{\{l:r\}} = 0$ for $r-l \leq 1$ and $F^*_l$ non-decreasing in $l,$
the minimum on the right hand side actually only has to be taken
over the values $l= 1, \ldots, r - 1.$
By the dynamic programming principle,
we successively compute  $F^*_1,$ $F^*_2,$ until we reach $F^*_N.$ 
As our primary interest are the optimal discontinuities $J,$ rather than the minimal functional values $F^*_N$,
we keep track of these locations.
An economic way to do so is to store at step $r$ the minimizing argument $l^*$ 
of \eqref{eq:recurrencePenalized} as the value $Z_r$ 
so that $Z$ encodes the boundaries of an optimal partition.
We refer to  \cite{friedrich2008complexity} for a more detailed description of the data structure and a visualization.
We record that the number of configurations that have to be checked for computing the minimum value is $O(N^2)$.

\subsection{Efficient computation of the spline energies for all intervals}\label{sec:efficient_energies}

We next develop an efficient procedure to compute the functional values 
\eqref{eq:epsilon_I} for all discrete intervals $\{ l:r\}$ with $1 \leq l \leq r \leq N.$
Before starting the development, we discuss why we need a specialized method for this specific task:
The algorithm of \cite{reinsch1967smoothing} aims for computing the minimizing argument of \eqref{eq:epsilon_I}, and the minimum value has to be computed from that minimizer. This procedure needs linear time for a single value of $\Ec_{\{l:r\}}$.
The method described by \cite{de1986efficient} gives direct access to the minimum energy $\Ec_{\{l:r\}}$ which is achieved by using orthogonal transformations. 
However, it is designed for fast computation of the energies for different $p$-parameters for signals of fixed length; no efficient algorithm for signals of increasing lengths is given. Hence, computing a single $\Ec_{\{l:r\}}$ needs linear time as well. Either approach would result in an algorithm of cubic complexity when used for computing all values $\Ec_{\{l:r\}}.$

We now develop a specialized method that 
updates the spline energy on the interval $\{ l:r\}$
to that on the interval $\{ l:r+1\}$ in constant time.
To keep the notation simple, we describe the procedure for $l = 1$ 
noting that the procedure works analogously for any other $l = 2, \ldots, N.$

Consider the functional in  \eqref{eq:epsilon_I} 
for an interval $I$ which contains exactly the data sites  $x_1, \ldots, x_r,$ where $2 \leq r \leq N-1,$
and denote its unique minimizing argument by $\hat f_I.$
The solution $\hat f_I$ is a cubic spline with natural boundary conditions which means that $p_i := \hat f_I|_{[x_i,x_{i+1}]}$ is a polynomial of degree at most $3$
for all $i=1, \ldots, r-1,$ that $\hat f_I$ is in $C^2(I),$ and that $\hat f_I''(x_1) = \hat f_I''(x_r) = 0$ \cite[Ch. XIV]{de1978practical}.
Let $f, f' \in \R^{r}$ be the vector of Hermite control points of $\hat f_I,$
meaning that  $f_i := \hat f_I(x_i)$ and, $f'_{i} := \hat f'_I(x_i),$ for $i=1, \ldots, r.$ 
Each polynomial $p_i$ is uniquely determined by $f_i, f'_{i}, f_{i+1},$ and $f'_{i+1};$ 
see e.g. \cite[p.61]{cheney1998introduction}. 
Hence $\hat f_I$ and its $2r$ Hermite control points $f, f'$ are one-to-one. 
We may express the polynomials $p_i$ in terms of the Hermite control points as 
\begin{equation}\label{eq:cubic_poly}
p_i(x) = \sum_{k=0}^3 c_k (x - x_i)^k	
\end{equation}
where $c_0 = f_i, c_1 = f'_i,$ 
	 $c_2 = -\frac{f'_{i+1} + 2f'_i}{d_i} + 3\frac{f_{i+1} - f_i}{d_i^2}, 
	c_3 = \frac{f'_{i+1} + f'_i}{d_i^2} + 2\frac{f_i - f_{i+1}}{d_i^3},$
and $d_i = x_{i+1} - x_i$ is the distance between two data sites \citep{dougherty1989nonnegativity}.
Using this representation, we obtain the integral of the square
of $p_i''$ on the interval $[x_i, x_{i+1}]$
in terms of $f_i, f'_{i}, f_{i+1},$ and $f'_{i+1}$: 
\[
\begin{split}
\int_{x_i}^{x_{i+1}} (p_i''(x))^2\, dx = & \frac{12 f_i^2}{d_i^3}-\frac{24 f_{i+1} f_i}{d_i^3}+\frac{12 f_{i+1}^2}{d_i^3} +\frac{12 f'_i f_i}{d_i^2}
+\frac{12 f'_{i+1} f_i}{d_i^2} \\ &  -\frac{12 f'_i f_{i+1}}{d_i^2}-\frac{12 f_{i+1} f'_{i+1}}{d_i^2}+\frac{4 (f'_i)^2}{d_i}+\frac{4 (f'_{i+1})^2}{d_i}+\frac{4 f'_i f'_{i+1}}{d_i}.	
\end{split}
\]
This is a quadratic form which can be written in matrix form as
$
	\int_{x_i}^{x_{i+1}} (p_i''(x))^2\, dx = v_i^T B v_i,
$
where 
$v_i = [f_i, f'_i, f_{i+1}, f'_{i+1}]$ and
\[
	B_i = \begin{bmatrix}
 \frac{12}{d_i^3} & \frac{6}{d_i^2} & -\frac{12}{d_i^3} & \frac{6}{d_i^2} \\
 \frac{6}{d_i^2} & \frac{4}{d_i} & -\frac{6}{d_i^2} & \frac{2}{d_i} \\
 -\frac{12}{d_i^3} & -\frac{6}{d_i^2} & \frac{12}{d_i^3} & -\frac{6}{d_i^2} \\
 \frac{6}{d_i^2} & \frac{2}{d_i} & -\frac{6}{d_i^2} & \frac{4}{d_i} \\
	\end{bmatrix}.
\]
The matrix $B_i$ is symmetric and positive semidefinite with two zero eigenvalues.
Invoking  a semidefinite variant of the Cholesky decomposition  we obtain the factorization
\begin{equation}\label{eq:Ai}
		 B_i = U_i^T U_i, \quad \text{with} \quad
U_i = \left[
\begin{array}{cccc}
 \frac{2 \sqrt{3}}{d_i^{3/2}} & \frac{\sqrt{3}}{\sqrt{d_i}} & -\frac{2 \sqrt{3}}{d_i^{3/2}} & \frac{\sqrt{3}}{\sqrt{d_i}} \\
 0 & \frac{1}{\sqrt{d_i}} & 0 & -\frac{1}{\sqrt{d_i}} \\
\end{array}
\right] \in \R^{2 \times 4}.
\end{equation}
For the next step, it is convenient to decompose $U_i$ as $ U_i = 
	[V_i, W_i ] $
	with $V_i, W_i \in \R^{2\times 2}$ and to define $e_1^T= [1,0] \in \R^{1\times 2}.$
	Let us define $A^{(r)} \in \R^{2r \times s},$ with $s = 3r -2,$ by \begin{equation}\label{eq:system_matrix}
A^{(r)} = \begin{bmatrix}
	\alpha_1 e_1^T & 0 & 0 & \cdots & 0\\
	 \beta V_1  & \beta  W_1 & 0 & \cdots & 0 \\
	0 & \alpha_2 e_1^T & 0 & \cdots & 0 \\
	0 & \beta V_2  & \beta W_2 &  \cdots & 0 \\
	0 & 0 & \alpha_3 e_1^T & \cdots & 0 \\
	\vdots & \ddots & \ddots & \ddots & \vdots  \\
	0 & \cdots  & 0 & \beta V_{r-1}  & \beta W_{r-1}  \\
	0 & \cdots  & 0  &	0 & \alpha_{r} e_1^T 
\end{bmatrix}
\end{equation}
with $\alpha_i = \frac{\sqrt{p}}{\delta_i} $ and $\beta = \sqrt{1-p}.$
Now the optimal functional value $\Ec_{\{1:r\}}$  can be expressed as follows: 
\begin{lemma}\label{lem:minimum_value}
Let $\Ec_{\{1:r\}}$ be given by \eqref{eq:epsilon_I} with $I ={\{1:r\}}.$ 
Then, 
\begin{equation}\label{eq:approx_error_lsq}
\Ec_{\{1:r\}} = \min_{u \in \R^{2r}} \|  A^{(r)} u - \tilde y^{(r)}\|_2^2,
\end{equation}
where $A^{(r)}$ is defined in \eqref{eq:system_matrix}, and $\tilde y^{(r)} \in 
\R^{s}$  is a vector of zeros except $\tilde y^{(r)}_{3i-2} = \alpha_i y_{i}$ for $i = 1, \ldots, r.$ 	
\end{lemma}
The proof is given in the supplementary material.

The form \eqref{eq:approx_error_lsq} has a key property for our purposes: it is a least squares problem in matrix form so that 
$A^{(r)}$ is a submatrix of $A^{(r+1)}$
and  $\tilde y^{(r)}$ is a subvector of $\tilde y^{(r+1)}.$
This submatrix relation allows us to update the QR-decomposition of the system $[A^{(r)}| \tilde y^{(r)}]$
to the QR-decomposition of the system $[A^{(r+1)}| \tilde y^{(r+1)}]$ by a constant number of Givens rotations.
Let $Q \in \R^{s \times s}$ orthogonal and $R \in \R^{s \times r},$ 
with $R_{1:r, 1:r}$ upper triangular and $R_{r+1:s, 1:r} = 0,$
such that
$A^{(r)} = QR,$ and let $z = Q^T \tilde y^{(r)}.$
Then the minimizer of \eqref{eq:approx_error_lsq}  is the solution of the first $r$ rows of the system 
\begin{equation}\label{eq:R_z}
	[R | z]
\end{equation}
and $\Ec_{\{1:r\}}$ is equal to the sum of squares of the last entries, $r+1, ..., \ldots, s,$ of the right hand side $Q^T \tilde y^{(r)}.$	
Passing from $r$ to $r+1,$ extends  the system by three rows and two columns:
\begin{equation}\label{eq:large_system}
	\left[
\begin{array}{ccc|c}
		 R_{1:s, 1:r-2}& R_{1:s, r-1:r}&0  & z \\
		 0 & \beta V_{r+1}  & \beta W_{r+1} &0 \\
		0  &	0 & \alpha_{r+1} e_1^T  & \alpha_{r+1} y_{r+1}
	\end{array}
	\right].
\end{equation}
This system can be brought to upper triangular form using 
Givens rotations.
By the upper triangular form of $R,$ 
the Givens rotations only act on the small subsystem
\begin{equation}\label{eq:small_subsystem}
		\left[
\begin{array}{cc|c}
         R_{r-1:r, r-1:r}&0 & z_{r-1:r} \\
		  \beta V_{r+1}  & \beta W_{r+1} &0 \\
			0 & \alpha_{r+1} e_1^T  & \alpha_{r+1} y_{r+1}
	\end{array}
	\right]
\end{equation}
where the left hand side is in $\R^{5 \times 4}.$
Having applied the Givens rotations, the system has the form
$
         [R'|  z']
$
with  $R' \in \R^{5\times 4}, $ $R'_{1:4,1:4}$ upper triangular and $R'_{5,1:4} =0.$
So the residual is coded in  the last entry of $z'.$ Hence we get the update
\begin{equation}\label{eq:update_eps_lr}
\Ec_{\{1:r+1\}} = \Ec_{\{1:r\}} + (z'_5)^2.	
\end{equation}
The update procedure is initialized (for the index $r = 2$)
by $\Ec_{\{1:2\}} = 0,$ and the system describing the initial state is given
$[R|Q^T \tilde y^{(2)}]$ where $QR = A^{(2)}$ is the QR-decomposition of $A^{(2)} \in \R^{4\times 4}.$

The above derivation leads to the following theorem.
\begin{theorem}\label{thm:update}
	The procedure described above (equations \eqref{eq:R_z} to \eqref{eq:update_eps_lr}) computes
	the array $[\Ec_{\{1:1\}}, \Ec_{\{1:2\}}, \ldots, \Ec_{\{1:N\}}]$ in $O(N)$ time complexity and $O(N)$ memory complexity.
\end{theorem}
The proof is given in the supplementary material.

For any other starting index $l$ with $1< l \leq r$
the values $\Ec_{\{l:r\}}$ are computed in the same fashion by calling the above procedure 
with the starting index $l.$

A note on the practical implementation seems useful at this point:
The described procedure uses Givens rotations for elimination.
Unfortunately, by their rowwise action, applying Givens rotations sequentially turns out to be relatively slow when implemented in standard linear algebra packages as used by Matlab. 
 However, as the Givens rotations applied to \eqref{eq:large_system}
only act on the subsystem \eqref{eq:small_subsystem},
we may use Householder reflections directly on \eqref{eq:small_subsystem}.
This allows to use  standard Householder-based QR decompositions.
Another possibility is a hardware-near implementation of the Givens rotation. For maintainability and readability reasons, we decided to use the standard linear algebra system as used by Matlab for the reference implementation.

If the sampling distances between adjacent data points $d_i$ have very different orders of magnitude, so that the global mesh ratio $({\max_{i=1,\ldots,N-1} d_i})/({  \min_{i=1,\ldots,N-1} d_i})$ is large,
the matrix $A^{(r)},$ defined by \eqref{eq:system_matrix}, may have a large condition number.
Then, solving \eqref{eq:approx_error_lsq} may get numerically unstable.
In such cases, one may try data binning
(meaning merging the data points with the smallest distances) to reduce the global mesh ratio.

\subsection{Compatibility with pruning strategies}\label{sec:pruning}

In the supplementary material, we show that the fast update strategy is compatible with the  specific order of the computation of two different pruning strategies, 
the PELT strategy from \cite{killick2012optimal} and the pruning from \cite{storath2014fast}, abbreviated here by FPVI.
Comparing the two strategies, the effectiveness of FPVI essentially grows  with the value of $\gamma$ whereas the effectiveness of PELT essentially grows with the number of detected discontinuities.

\subsection{Overall algorithm and computational effort}\label{sec:overall}

Putting together the procedures described in the last subsections, 
we obtain an algorithm which computes a solution  $\hat J$ of the minimization problem \eqref{eq:pcw_spline}; that is,  $\hat J$ is a global minimizer  of the target function in \eqref{eq:pcw_spline}.
Having computed the optimal discontinuity set $\hat{J},$ the corresponding CSSD $\hat{f}$ 
is computed using the algorithm of \cite{reinsch1967smoothing}.
This completes the solution of \eqref{eq:pcw_spline_lagrange}.
Altogether, we obtain:
\begin{theorem}\label{thm:algorithm}
	The algorithm described above (sections \ref{sec:reduction} to \ref{sec:overall})
	computes a solution of \eqref{eq:pcw_spline_lagrange}, i.e., a global minimizer of the target function in \eqref{eq:pcw_spline_lagrange}.
	The worst case time complexity is $\Oc(N^2)$ and the memory complexity is $\Oc(N).$
\end{theorem}
The proof is given in the supplementary material.

For details on the implementation, we refer to the provided reference implementation, cf. Section~\ref{sec:implementation}.

 We point out that the scaling in runtime is  often more favourable than the worst case scenario. This is an effect of the pruning strategies described in section \ref{sec:pruning}.
We observed in our experiments that the runtime scales approximately linearly
if the number of detected discontinuities grows linearly with the signal's length. 
If there are no or very few detected discontinuities the runtime grows approximately quadratically.
We refer to the experimental section for an illustrating example.
 Further, \cite{killick2012optimal} made similar observations for piecewise estimators and gave theoretical justifications for the PELT pruning strategy.

\begin{remark}
It is straightforward to extend the method to work for vector-valued data
meaning that  $y_i \in \R^D$  for all $i = 1, \ldots, N.$
In this case, $\Ec_{\{l:r\}}$ is simply the sum of the energies over all vector components $\sum_{j=1}^D \Ec_{\{l:r\}, j},$
where $\Ec_{\{l:r\}, j}$ is the spline energy of the $j$-th component of $y$ on $\{l:r\}.$
It is clear that the computational effort scales linearly with the dimension~$D.$	
\end{remark}

\subsection{Uniqueness result}\label{sec:uniqueness}
As discussed in Section~\ref{sec:reduction}, the precise location of discontinuities between data sites does not influence the penalty function, so that a solution can only be unique up to shifts of discontinuities in between the data sites.
Having reduced the potential discontinuities to midpoints (cf. Equation \eqref{eq:pcw_spline_discrete}), the optimal discontinuity set $\hat{J}$ may still not be unique.  A simple example is the data $x = [0,1,2]$ and  $y = [0,1,0].$
For sufficiently small $\gamma > 0,$ both $J = \{\frac{1}{2}\}$ and $\tilde J = \{\frac{3}{2}\}$ 
are optimal because they result in the same target functional value $\gamma.$
A similar situation is found whenever the signal has odd length
and has $\ceil{N/2} - 1$ discontinuities. 
As it can be seen from the  proof of Lemma~\ref{lem:J_star}, 
a solution with $\ceil{N/2} - 1$ discontinuities
is a piecewise linear function interpolating at most two data points.
The ambiguities in the examples can be attributed to segments of length two or one,
on which a smoothing spline acts  interpolatory.
If partitions with  segments of length less than $3$ can be ruled out, we obtain a uniqueness result for a.e. input $y \in \R^N$ (in a Lebesgue almost everywhere sense).
The result is inspired by a related theorem of \cite{wittich2008complexity} on the piecewise constant model.

\begin{theorem}\label{thm:uniqueness}
	We consider the minimization problem \eqref{eq:pcw_spline_lagrange} for the CSSD in its formulation 	
	\eqref{eq:pcw_spline_discrete} (with discontinuity set reduced to midpoints).
	Then,  
	the minimizing values $f(x_i)$ at the data sites $x_i$ are uniquely determined  for a.e. input $y \in \R^N$ (w.r.t. Lebesgue measure.)
	Further, minimizing segments $I^\ast$ of length at least three are uniquely determined for (Lebesgue-) a.e. input $y \in \R^N.$
	In particular, the minimizing function $f$ is unique in $ [ x_{\min}, x_{\max} ]$ 
	where $x_{\min}, x_{\max}$ 
	denote the minimal and maximal argument value of the $x_i$ in the minimizing segment $I^\ast$ of length at least three.
	Consequently, minimizing partitions with all segments being of length at least three and corresponding minimizing functions (outside the intervals with the jumps) are uniquely determined  for  a.e. input $y \in \R^N.$
\end{theorem}
The proof is given in the supplementary material.

In the present setup, a minimizer of \eqref{eq:pcw_spline_discrete} with at least 3 data sites between two discontinuities is obtained if  $\gamma $ is chosen sufficiently large. 
We note that 
one could also constrain the search space of \eqref{eq:pcw_spline_discrete} w.r.t. the size of the discontinuity sets
as an alternative approach.

\section{Experiments and parameter selection}

\subsection{Implementation and setup}\label{sec:implementation}
The proposed algorithm was implemented in Matlab R2022b. 
For implementation details,
we refer to the commented source code provided on Github\footnote{Reference implementation of CSSD provided at \url{https://github.com/mstorath/CSSD}.}.
The numerical experiments were performed on a laptop with 2.4 GHz 8-Core Intel Core i9, 32 GB RAM.
The choice of the model parameters $p$ and $\gamma$ depends on the sampling distances on the abscissa,
on the scale of the data on the ordinate, and, of course, on the shape of the underlying signal.
We set $\delta_i = \sigma$  for the synthetic experiments, where $\sigma$ is the standard deviation of the noise, 
and $\delta_i = 1$ for the other experiments, for $i= 1, \dots, N.$

\subsection{Simulated experiments with a manually selected parameter set.}
As the proposed method is designed for estimating smoothly varying functions with discontinuities,  we primarily choose  test signals with these properties.
A classical signal of this class is the \enquote{HeaviSine} function used in the work of  \cite{donoho1994ideal}. 
It is given by  $g_2(x) :=  4\sin(4\pi x) - \sign(x - 0.3) - \sign(0.72 - x).$ 
As second function of this class, 
we use the function $g_1(x) := J_1(20x) + x \cdot {1}_{[0.3, 0.4] }(x) - x \cdot {1}_{[0.6, 1] }(x)$ as test signal, where $J_1$ is the Bessel function of the first kind and ${1}_{S}$ is the indicator function on $S$. Compared to the HeaviSine, it has one discontinuity more and a varying amplitude. The discontinuity locations of $g_1$ are chosen such that $g_1$ shares  one discontinuity with $g_2$ and such that the segments have different lengths.
For the experiment shown in Figure~\ref{fig:synthetic} (in the introduction),
we sampled the signal at $N = 100$ random points  in the interval $[0,1]$ (uniformly distributed)  and corrupted them by additive zero mean  Gaussian noise of standard deviation $\sigma = 0.1.$
For a classical spline, a visually reasonable result was obtained using the parameter $p=0.999.$ 
Using this $p$-parameter, we report the results for several $\gamma$-values.

Experimental results for vector-valued input data are given in the supplementary material.

As the resulting function is a cubic spline with discontinuities
	it can be described by a piecewise cubic polynomial with discontinuities in $J.$
	Therefore, the result may serve as functional description of the discontinuous signal.

\subsection{Computational costs}\label{sec:comp_costs}
 Next, we investigate the computational costs of the algorithm.
 As mentioned in the introduction,  there is no  freely available software package for solving  \eqref{eq:pcw_spline_lagrange} or its reduced form  \eqref{eq:pcw_spline_discrete}, to the authors knowledge. 
To assess the benefits of the proposed method, we implemented a baseline solver which relies on standard libraries for solving the reduced form \eqref{eq:pcw_spline_discrete}. To this end, we combine a dynamic programming approach  with the PELT-pruning strategy \citep{winkler2002smoothers, killick2012optimal}, implemented in the Python module \texttt{ruptures}\footnote{\url{https://github.com/deepcharles/ruptures}} \citep{truong2020selective},
for computing the optimal partition in \eqref{eq:pcw_spline_discrete}. The necessary values of $\Ec_I$  on each interval (cf. Equation \eqref{eq:epsilon_I}) were obtained as follows. The smoothing spline module \texttt{csaps}\footnote{\url{https://github.com/espdev/csaps}}, which is described as a Python replication of De Boor's corresponding Fortran method, was used to compute a minimizer $\hat f_I$ of \eqref{eq:epsilon_I} 
	which in turn was used to compute $\Ec_I$ using the corresponding equation given by \cite{de1978practical}.
	As this baseline solver does not use the fast update scheme developed in this work,
	its worst case time complexity  is $O(N^3)$ because the number of discrete intervals grows quadratically and the computational effort for a smoothing spline grows linearly in the interval length.

For the runtime comparison, we consider again the HeaviSine signal $g_2$ used  in the last experiment,
and we distinguish two experimental setups:
In the first scenario, we increase the number of data points $N$ 
by increasing the number of samples of the HeaviSine 
so that the number of discontinuities remains two.
In the second scenario, we increase the number of data points
by repeating the HeaviSine signal so that the number of discontinuities grows linearly with the signal length.
We computed the CSSD for the parameters $p=0.9999$ and $\gamma = 20$
which resulted in reasonably good estimates for the underlying true signal.

\begin{table}[t] 
\centering
		\begin{tabular}{lrrrrrrr} 
\toprule  
Signal length $N$ & 250 & 500 & 1000 & 2000 & 4000 & 8000\\ 
 \toprule 
\multicolumn{7}{l}{\textbf{Constant number of discontinuities}}\\ 
 \midrule 
Baseline & $1.1 \cdot 10^{6}$ & $2.7 \cdot 10^{6}$ & $2.0 \cdot 10^{7}$ & $1.5 \cdot 10^{8}$ & $1.3 \cdot 10^{9}$ & $7.3 \cdot 10^{9}$\\ 
  \cmidrule{2-7} 
&45.7 & 110.1 & 453.4 & 1909.1 & 8704.0 & 33331.3\\ 
  \cmidrule{1-7} 
Prop. + FPVI & $2.5 \cdot 10^{4}$ & $10.0 \cdot 10^{4}$ & $4.0 \cdot 10^{5}$ & $1.7 \cdot 10^{6}$ & $7.0 \cdot 10^{6}$ & $2.9 \cdot 10^{7}$\\ 
  \cmidrule{2-7} 
&0.4 & 1.0 & 4.2 & 17.3 & 72.6 & 297.1\\ 
  \cmidrule{1-7} 
Prop. + PELT & $1.8 \cdot 10^{4}$ & $4.3 \cdot 10^{4}$ & $1.7 \cdot 10^{5}$ & $6.8 \cdot 10^{5}$ & $2.7 \cdot 10^{6}$ & $9.5 \cdot 10^{6}$\\ 
  \cmidrule{2-7} 
&0.3 & 0.6 & 2.2 & 8.7 & 34.5 & 120.9\\ 
 \midrule 
 \multicolumn{7}{l}{\textbf{Increasing number of discontinuities}}\\ 
 \midrule 
Baseline &$3.9 \cdot 10^{5}$ & $8.0 \cdot 10^{5}$ & $7.4 \cdot 10^{5}$ & $1.4 \cdot 10^{6}$ & $2.9 \cdot 10^{6}$ & $5.7 \cdot 10^{6}$\\ 
  \cmidrule{2-7} 
&28.6 & 57.8 & 81.8 & 158.9 & 323.6 & 652.3\\ 
  \cmidrule{1-7} 
Prop. + FPVI &$2.1 \cdot 10^{4}$ & $9.1 \cdot 10^{4}$ & $1.6 \cdot 10^{5}$ & $2.4 \cdot 10^{5}$ & $5.0 \cdot 10^{5}$ & $1.3 \cdot 10^{6}$\\ 
  \cmidrule{2-7} 
&0.2 & 1.0 & 1.7 & 2.5 & 5.3 & 13.7\\ 
  \cmidrule{1-7} 
Prop. + PELT &$1.2 \cdot 10^{4}$ & $2.4 \cdot 10^{4}$ & $3.5 \cdot 10^{4}$ & $6.9 \cdot 10^{4}$ & $1.4 \cdot 10^{5}$ & $2.8 \cdot 10^{5}$\\ 
  \cmidrule{2-7} 
&0.2 & 0.3 & 0.5 & 1.0 & 2.0 & 3.9\\ 
 \bottomrule 
  	\end{tabular}	
	\caption{Computational effort of the baseline solver and the proposed solver with different pruning strategies for computing a CSSD.
	The top rows represent the frequency a data item is visited by the algorithm,
	and the bottom rows the total runtime in seconds, respectively.
	 The subtables show two scenarios:
	 Increasing sampling density which results in a constant number of discontinuities,
	 and repeating the signal so that the number of discontinuities scales linearly with the signal length. 
	 The reported values are averages over three runs with different noise realizations
	 and a fixed parameter set.
	}
		\label{tab:runtimes}
\end{table}

In Table~\ref{tab:runtimes}, we report two measures of the computational effort,
the runtimes and the frequency 
each data point is \enquote{visited} for computing all necessary $\Ec_{[l:r]}.$
(The baseline algorithm needs the vector $y_{l:r}$ of length $r-l+1$
while the proposed method only needs a single data point, $y_l$ or $y_r,$ depending on the order of update.) 

We start our analysis by examining the frequencies at which each data point is visited. In the first scenario, where the number of discontinuities is held constant, the baseline method shows a  pattern more akin to cubic growth, whereas the proposed method rather shows quadratic growth. For the second scenario, the growth in $N$ trends towards linear for all considered methods. However, the hidden constant depends on the average segment length, and the relation is quadratic for the baseline method, and linear for the proposed method.
 Importantly, in both scenarios, the proposed method significantly reduces computational effort with respect to the total frequencies at which each data point is visited.

We next look at the corresponding execution times. 
In the first scenario, we note that the corresponding runtimes 
for the baseline method exhibit a more favorable growth pattern 
compared to the growth pattern of the frequencies. 
This appears to be a consequence of the fact that
the \texttt{csaps} function, which is used for computing the smoothing splines, shows a predominantly constant growth, and this diverges from the asymptotic linear growth for the majority of the signal lengths evaluated by the algorithm. An asymptotic cubic growth in runtime is expected for even longer signals than the ones currently used for comparison.
Overall, the proposed method is much faster 
than the baseline method.
In particular, the proposed methods allows to  process signals of moderate size on 
a laptop repeatedly within a reasonable time frame.

In both scenarios of the above example,  the PELT pruning is more effective 
than the FPVI pruning. Similar observations were made for other signals (not shown in this paper) in  which the CSSD detects at least a few discontinuities. 
Interestingly, for the parameter selection by cross validation on the real data sets (Section~\ref{sec:real_data}) it is the other way round: the runtimes with FPVI pruning were 2.1 min (geyser data) and 5.8 h (stock data) 
whereas runtimes with PELT pruning were 2.7 min  and 10.7 h, respectively.
(A detailed table of the runtimes for each sampled $(p, \gamma)$ value of the Geyser data is given in the supplementary material.)
The reasons for the lower runtimes of FPVI pruning in these cases are that a relatively low number of discontinuities are detected  and that the optimization procedure frequently evaluates the target functional for large values of $\gamma$  where FPVI is more effective than PELT.

\subsection{Automatic parameter selection} 
As for classical smoothing splines, it may be sufficient for some practical purposes to choose the 
model parameters  by visual inspection;
yet a procedure for automatic parameter selection is useful (cf. the discussion on this topic by \citet[Sec.~4]{silverman1985some}).
To obtain an automatic proposal for the model parameters  $p$ and $\gamma,$
we here use K-fold cross validation (CV) as follows.
We partition the data randomly into $K$ folds of approximately equal size.
The K-fold cross validation score for $p$ and $\gamma$ is given by
$
	\mathrm{CV}(p, \gamma) = \frac{1}{N}\sum_{k=1}^K  \sum_{i \in Fold_k} ((\hat f^{-k}_{p, \gamma}(x_i) - y_i)/\delta_i)^2.
$
Here,  $\hat f^{-k}_{p, \gamma}$
denotes the result of the proposed algorithm applied to the all data  except those in fold $k.$ 
(Recall from Section~\ref{sec:uniqueness} that the solutions of \eqref{eq:pcw_spline_lagrange} might be non-unique. Yet, the proposed algorithm provides a unique solution, namely the one with the largest possible right-most interval, followed by the largest possible penultimate interval, and so on.)
If an evaluation point $x_i$ coincides with a discontinuity of $f^{-k}_{p, \gamma},$ we take the mean of the lefthand and the righthand limits: $f^{-k}_{p, \gamma}(x_i) = \lim_{t \to 0+}\frac{1}{2}(f^{-k}_{p, \gamma}(x_i-t) + f^{-k}_{p, \gamma}(x_i+t)).$
A typical choice is $K = 5$ folds \citep{hastie2009elements} which we adopt here.
Optimal parameters $p, \gamma$ in the sense of K-fold cross validation 
are minima of the scoring function $(p, \gamma) \mapsto \mathrm{CV}(p, \gamma).$
To find a good parameter set in the sense of CV, we improve a starting value $(p_0, \gamma_0)$ using  standard derivative-free optimizers: a global search via simulated annealing  followed by a local refinement using the Nelder-Mead  method. We use the implementations of Matlab with default options. 
The simulated annealing algorithm utilizes a balance of exploration (looking for new, possibly better solutions) and exploitation (optimizing around the best solution found so far), regulated by the temperature and its cooling schedule. By the mechanism of being able to escape local minima, it is less sensitive to the starting value than local methods. To save computation time, it it reasonable to use a starting value that gives a visually reasonable result based on the domain knowledge or prior experience.
To enhance parallel processing,  simulated annealing could be replaced by a grid search of the $(p, \gamma)$ domain.
 For the optimization, $\gamma$ is parametrized by  $\gamma = p \frac{q}{1-q}$ 
with $q \in [0,1).$ (This parametrization is obtained by dividing the functional by $p$ and reparametrizing the resulting penalty $\gamma' = \gamma/p$ by $\gamma' = \frac{q}{1-q}$).

	The procedure gives reasonable results in the conducted experiments (see results further below), but it comes with the limitations that the standard optimizers often need a high number of function evaluations leading to long runtimes,
and that the resulting parameters are not guaranteed to be global minimizers of the CV scoring function.
The second limitation can be mitigated by restarting the optimization with different starting values until no further improvement of the CV score is observed. 
As the costs for this improvement are additional computation time and/or manual refinement,
we use this strategy only for the two real data experiments presented further below.
For the  synthetic example described further below, the initial values 
$p_0 = 0.99$ and $\gamma_0 =1$ are used. (These parameters were chosen based on their ability to produce visually acceptable results for a set of sample signals.)

We compare the results of the proposed method to the following baseline method. We use the Bayesian ensemble method of \cite{zhao2019detecting}, implemented in the toolbox \texttt{beast}\footnote{Source code retrieved from \url{https://github.com/zhaokg/Rbeast} on Nov. 2, 2022.}, to estimate changepoints and a corresponding piecewise estimate of the signal. To this end, the function \texttt{beast\_irreg} (variant for nonuniform sample distances) is called with deactivated seasonality option (because the signal has no seasonality), time step option $0.01$ (binning the non-equidistant samples to 100 bins of same length), and default parameters otherwise. 
	  (As internal preprocessing, \texttt{beast\_irreg}  aggregates data over bins of fixed width to obtain equidistant data sites. The procedure returns the first bin locations after the corresponding change; so to match the \enquote{midpoint convention} used in this paper, we work with the midpoint between the returned bin and the preceding non-empty bin.)
We then fit a piecewise interpolating spline to the estimated values between the changepoints to obtain a functional estimate on the entire domain $[0,1]$.
	(The interpolating spline is preferred over the smoothing spline here, because a smoothing spline would give different $y$-values at the data sites than those estimated by \texttt{beast} and would introduce an additional parameter.)
	
	\begin{figure}[]
	\includegraphics[width=\textwidth]{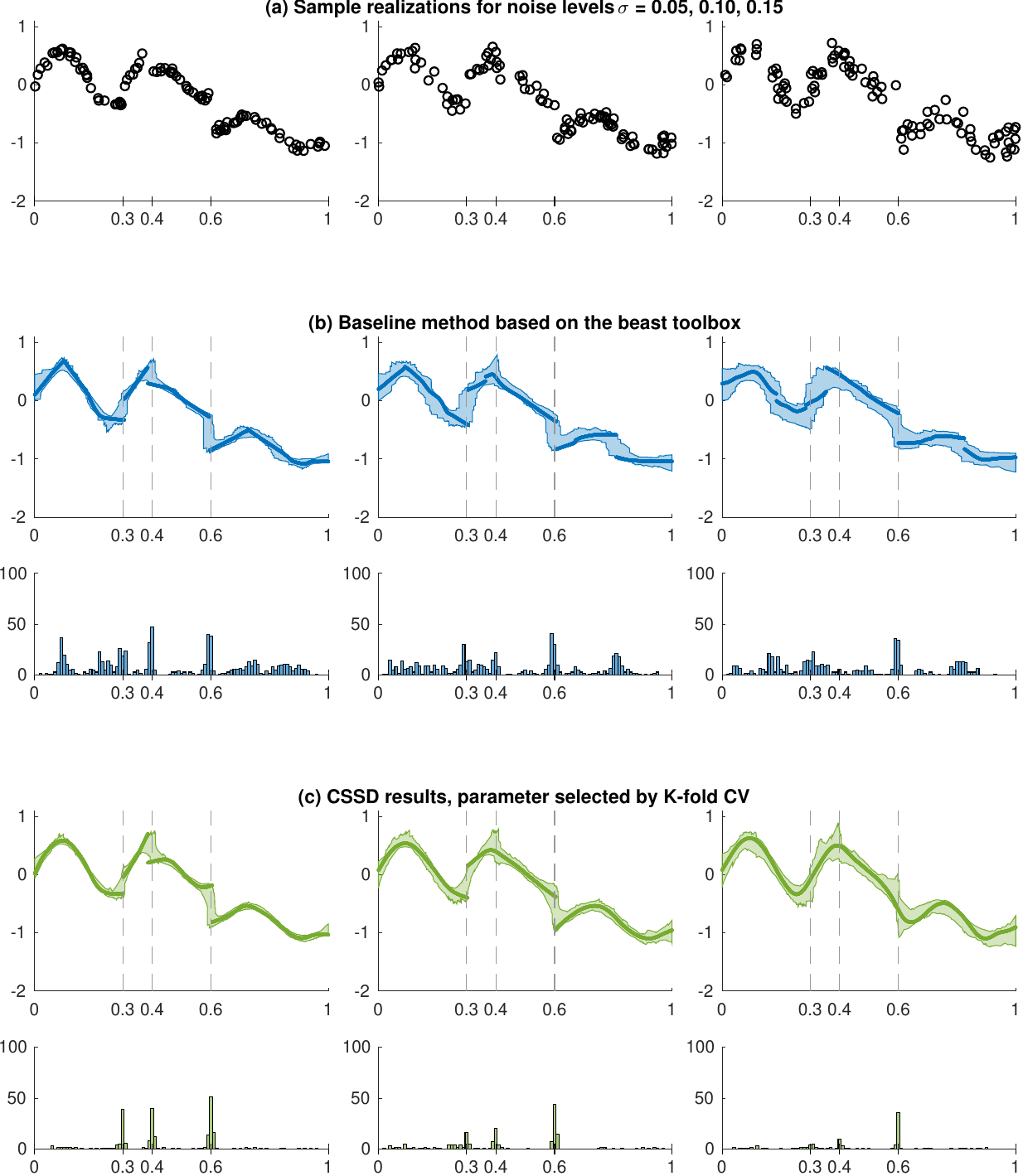}
	\caption{CSSD results with automatically selected parameter based on K-fold CV in comparison with a baseline method based on the beast toolbox. 
	\textit{(a)} Noisy samples of the test signal $g_1$ shown in Figure~\ref{fig:synthetic}a with the noise levels $\sigma = 0.05, 0.1, 0.15$ (from left to right). \textit{(b)} Results of the baseline method consisting of a method for piecewise regression with automatic parameter selection and subsequent spline interpolation. The dashed vertical lines indicate the true discontinuity locations. (As in Figure~\ref{fig:synthetic},
	the solid line represents the result for the sample signal, and the shaded regions represents the $2.5 \%$ to $97.5 \%$ quantiles, and the histograms represent the frequencies of the discontinuity locations.)
	\textit{(c)} Results of the CSSD with automatic parameter selection based on K-fold cross validation. 
}
	\label{fig:synthetic_CV_Equi}
\end{figure}

	We use the experimental setup used for Figure~\ref{fig:synthetic} and two additional noise levels to compare the results of the two methods.
	Figure~\ref{fig:synthetic_CV_Equi} shows the results for $100$ signal realizations.
	We observe that the baseline method returns principally more changepoint candidates than the proposed method, and that a considerable portion of them are located away from the true discontinuity locations, e.g. near the local extrema of the test signal. The proposed method exhibits less spurious discontinuities than the baseline method and the cluster points in the histograms are more sharply located around the true discontinuity locations.

\subsection{Results on real data using automatic parameter selection}
\label{sec:real_data}
\begin{figure}[t]
	\centering
	\includegraphics[width=0.9\textwidth]{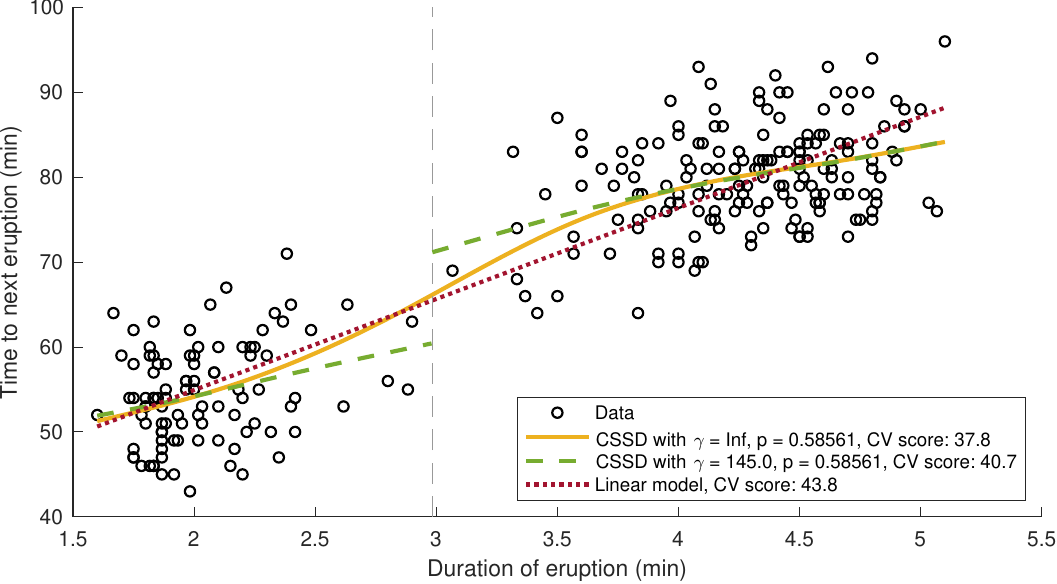} 
	\caption{	Fitting a CSSD to the Old Faithful data (circles):
	If the parameter is selected based on K-fold CV
	we obtain a result without discontinuities which coincides with a classical smoothing spline (solid curve).
	Keeping the selected $p$-parameter and lowering the $\gamma$ parameter sufficiently gives a two-phase regression curve (dashed curves) with a breakpoint near $x = 3$ (dashed vertical line), and the two curve segments are nearly linear. 
	Both of the above parameter sets yield better CV-scores  than a linear model (dotted line).
	} 
	\label{fig:geyser}
\end{figure}

The next experiment follows an example presented in the work of \cite{silverman1985some} on smoothing splines.
The data set contains the duration of eruptions along with the waiting time to the next eruption of the Old Faithful geyser in the Yellowstone National Park, USA; 	see Figure~\ref{fig:geyser}.\footnote{The data set was retrieved from the supplementary material of \cite{wasserman2004all}  \url{https://www.stat.cmu.edu/~larry/all-of-statistics/=data/faithful.dat}}
When fitting a single straight line to the data, curvature effects in the residuals are observed. \cite{silverman1985some} argues that fitting a smoothing spline is a useful exploratory step towards the choice of a reasonable model. Inspecting the shape of the resulting spline suggest a two phase linear regression model as plausible alternative.
	Remarkably, the CSSD model comprises both of the above models:
	The automatic parameter selection results in $\gamma = \infty$  so that the corresponding CSSD coincides with a classical smoothing spline.
	An appropriately chosen smaller $\gamma$-parameter gives a two phase model which is nearly linear on the two segments, and tends to a piecewise linear model for $p \to 0.$ 
	Ranking these three models by their CV scores, the classical spline can be considered as the best one, the linear model as the worst, and the two phase model lies in the middle.
	The example illustrates that the CSSD can serve as a useful tool for exploratory data analysis.

Eventually we report the results on a stock market time series. 
The present data shows the logarithm of the closing price of the Meta/Facebook stock 
from May 18, 2012, 
	to May 19, 2022.\footnote{Data source: \url{https://www.macrotrends.net/stocks/charts/FB/meta-platforms/stock-price-history}}
	The dates of the discontinuities can be related to strong market reactions after business reports:
	For example, on July 24, 2013, Facebook announced remarkable rises in revenues \citep{wilhelm2013facebook}. On July 25, 2018, Meta announced  lower profit margins \citep{reuters2018meta},
 and on February 2, 2022, the price of Meta shares dropped strongly 
	as Facebook reported a decline in daily users for the first time \citep{reuters2022meta}.
	In this scenario, detected discontinuities could be interpreted as abrupt changepoints of the signal. The computation times are reported at the end of Section~\ref{sec:comp_costs}.
\begin{figure}[t]
	\centering
	\tikzstyle{myspy}=[spy using outlines={black,lens={scale=3},width=0.6\textwidth, height=0.27\textwidth, connect spies, thick,
every spy on node/.append style={very thick}}]
	
	\def\figwidth{0.48\columnwidth}
\def\nodeSpy{(1.6, 2.55)}
\def\nodeWindow{(6.7,-0.5)}
\def\nodeSpyB{(1.7, 2.)}
\def\nodeWindowB{(2.2,-0.5)}
  	\begin{tikzpicture}[myspy]
	\node {\includegraphics[width=1\textwidth]{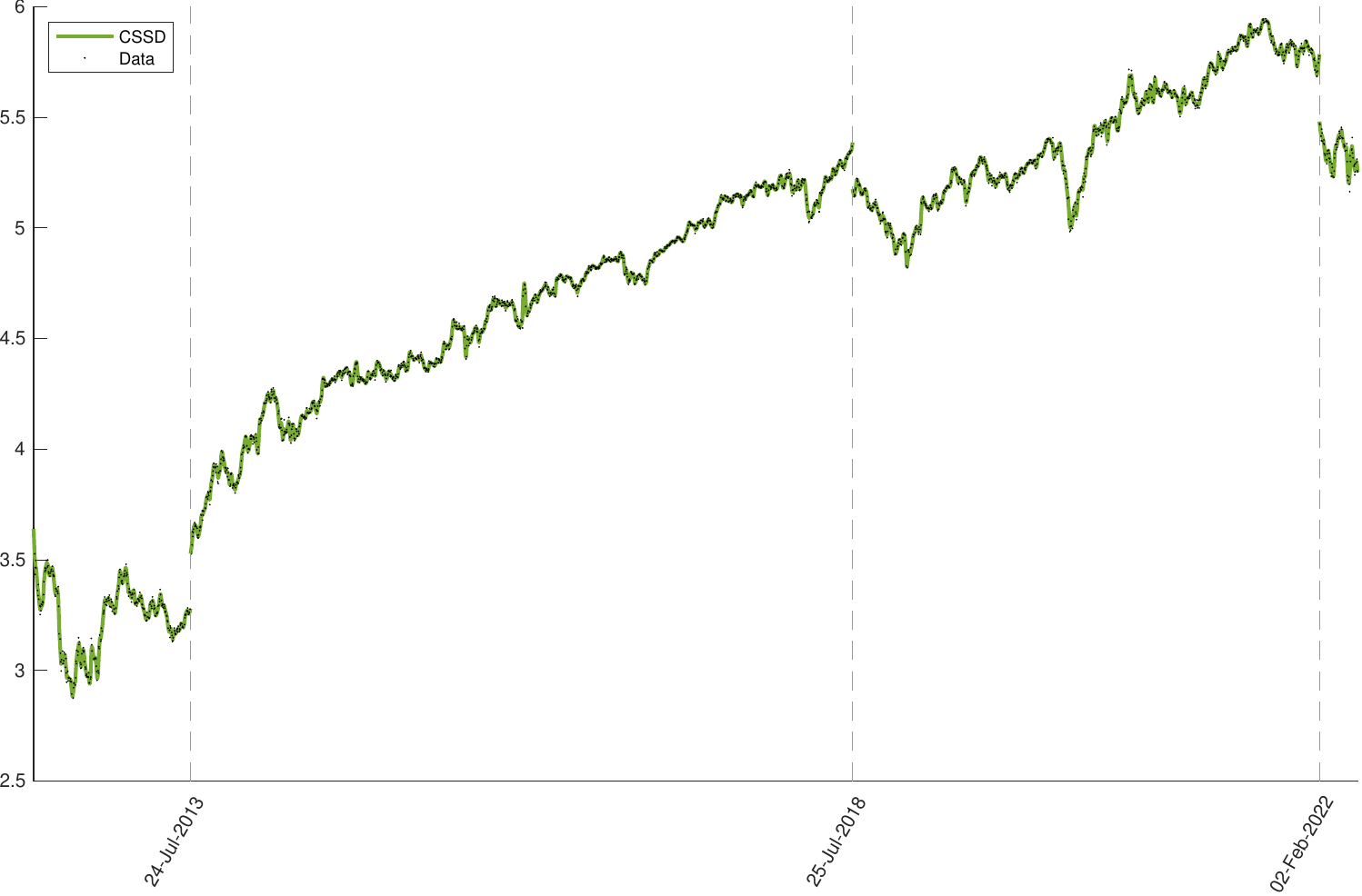}};
	\spy on \nodeSpy in node [left] at \nodeWindow;
\end{tikzpicture}	
	\caption{The dots represent the logarithm of the closing prices of the Meta stock from May 18, 2012, 
	to May 19, 2022. The curve represents the CSSD with parameters determined by K-fold CV ($p = 0.4702,$ $\gamma = 0.0069$). The dashed vertical lines indicate the discontinuities of the CSSD, and the ticks correspond to the date before the discontinuity.}
	\label{fig:stock}
\end{figure}

\section{Conclusion and outlook}
We have studied a variational model for cubic smoothing splines with unknown discontinuity locations which is a special case of the weak rod model.
This model comprises classical continuous cubic splines (for sufficiently large $\gamma$ parameters) and piecewise linear regression (for $p \to 0$). 
We have shown that the solutions are unique in an almost everywhere sense on segments comprising at least three data sites.
We developed an efficient algorithm for computing 
the global minimizer of the model's underlying optimization problem.
The runtime experiments indicated that the method is particularly efficient for signals with many discontinuities. 
But also in the general case, the algorithm can be executed on a laptop in reasonable time for signals of moderate lengths.
The algorithm is applicable to data with non-equidistant sampling points and to vector-valued data.
Automatic parameter selection based on K-fold cross validation gave reasonable results for synthetic and real data; the main limitation of this selection strategy is currently the relatively high runtime.
The numerical examples have illustrated potential applications of a CSSD:
 as function estimator for a discontinuous signal, as basis for a changepoint detector, or
as tool for exploratory data analysis.  

Open questions for future work include faster algorithms for automatically selecting the model parameters,
extension to splines of higher orders and higher dimensions, and further investigation of theoretical properties of the estimator.

\section*{Supplementary Material}
\begin{description}
      \item[Supplementary document:] The supplementary document contains proofs of the theorems, justifications for the compatibility with the pruning strategies, and additional experiments.
 \end{description}

\section*{Acknowledgement}

MS  was supported by the research program \enquote{Informations- und Kommunikationstechnik} of the Bavarian State Ministry of Economic Affairs, Regional Development and Energy  (DIK-2105-0044 / DIK0264). 
AW acknowledges support of Deutsche Forschungsgemeinschaft (DFG) under project number 514177753.
The authors report there are no competing interests to declare.

{\footnotesize
\bibliographystyle{myplainnat}
\bibliography{pcwSpline}
}

\newpage
\setcounter{page}{1} 

\appendix

\section{Supplementary document for \enquote{Smoothing splines for discontinuous signals}}

This is supplementary material for the paper: Martin Storath, Andreas Weinmann, \enquote{Smoothing splines for discontinuous signals}.
The  equations, statements and references refer to the corresponding 
elements in the main document. 

\subsection{Proofs}

\subsubsection{Proof of Lemma \ref{lem:J_star}}\label{app:lemma1}

\begin{proof}
1. Let $J$ be a minimizer of \eqref{eq:pcw_spline_lagrange}.
Assume that $J$ has two elements $j_1, j_2,$ with $j_1 <j_2,$ between two adjacent data sites $x_i, x_{i+1}$.
Let $f_J$ be the corresponding optimal function.
By a basic property of cubic smoothing splines (see e.g. \citep{silverman1985some}),
 $f_J$ is an affine linear function on $[x_i, j_1).$ This affine linear function can be extended to $[j_1,j_2).$
The extension does not increase the roughness penalty, and decreases the discontiunity penalty by $\gamma > 0.$ Hence,  $J$ was not a minimizer.
For more than two discontinuities, the above argument is applied repeatedly.
	
	2. For data of length $N,$ we find a discontinuity set $J$ of size $\ceil{N/2} - 1$
	which partitions  $[x_1, x_N]$  into  $\ceil{N/2}$ segments containing at most two data sites each.
	On each of the segments, the optimal spline is a linear function connecting the two data points
	(or a constant function if the segment contains only one point)
	so that both the data penalty and the roughness penalty become $0.$
	Thus, the only remaining penalty is the discontinuity penalty which is equal to
	 $\gamma J = \gamma (\ceil{N/2} - 1).$	
	 This shows that the minimum functional value of \eqref{eq:pcw_spline_lagrange} is bounded from above by $\gamma (\ceil{N/2} - 1).$ Now if $|\hat{J}| > (\ceil{N/2} - 1),$ the functional value is strictly larger than $\gamma (\ceil{N/2} - 1)$ which is a contradiction.
\end{proof}

\subsubsection{Proof of Lemma \ref{lem:minimum_value}}\label{app:lem2}

\begin{proof}
Denote the functional in \eqref{eq:epsilon_I}, with $I = \{1:r\},$
by $L_1(f)$ (defined for arguments $f \in C^2([x_1,x_r])$),  and denote the functional in \eqref{eq:approx_error_lsq}
by $L_2(u)$ (defined for arguments $u \in \R^{2r}$).

Let $\hat f$ be the minimizer of $L_1,$
and let $\hat v \in \R^{2r}$ such that $\hat v_{2i -1} = \hat f(x_i)$  and $\hat v_{2i} = \hat f'(x_i)$ for all $i=1,\ldots,r.$
In the derivation before Lemma \ref{lem:minimum_value} (equations \eqref{eq:cubic_poly} to \eqref{eq:system_matrix})
we constructed $A^{(r)}$ and ${\tilde y}^{(r)}$
such that $\Ec_{\{1:r\}} = L_2(\hat v).$
This immplies that
\begin{equation}\label{eq:lem_eps_1}
	 \Ec_{\{1:r\}} = L_2(\hat v) \geq  \min_{u \in \R^{2r}}L_2(u) = L_2(\hat u).
\end{equation}
	where $\hat u \in \R^{2r}$ is the minimizer of $L_2.$

	 Let $H$ be the unique piecewise cubic polynomial described by $\hat u,$
	i.e. $H(x_i) = \hat u_{2i -1}$  and $H'(x_i) = \hat u_{2i}.$
	By construction $H$ is a piecewise cubic polynomial and $ H \in C^1([x_1,x_r]).$
	 Let $S$ be the interpolating cubic spline of the points $(x_i, H(x_i))$ for $i=1,\ldots, r$ with natural boundary conditions.
	 Since $H(x_i) = S(x_i)$ for all $i=1,\ldots, r,$ $H$ and $S$ have the same data deviation penalty
	 (the first term in \eqref{eq:epsilon_I}).
	 	Towards a contradiction, assume that $H \neq S$ on $[x_1,x_r].$
	 	 $H'$ is absolutely continuous, because it is continuous and differentiable up to a finite set.
	Thus $H'$ is differentiable in the sense of Lebesgue, and  $H'' \in L^2([x_1,x_r]),$ because $H''$ is piecewise continuous with only step discontinuities at the data sites. Because $H \neq S$ by assumption, we get by the second minimality property of splines \citep{schoenberg1964best} that
	\[
		\int_{x_1}^{x_r}	 (H''(x))^2 \, dx > \int_{x_1}^{x_r}	(S''(x))^2 \, dx.
	\]
	Therefore, we get for $v \in \R^{2r}$ with $v_{2i -1} = S(x_i)$  and $v_{2i} = S'(x_i)$ that
	$
		L_2(v) < L_2(\hat u)
	$
	which contradicts the optimality of $\hat u.$ Hence, $H = S.$ 
	It follows that 
	\begin{equation}\label{eq:lem_eps_2}
		L_2(\hat u) = L_1(H) = L_1(S) \geq \min_{f \in C^2([x_1, x_r])} L_1(f) = \Ec_{\{1:r\}}.		
	\end{equation}
	Combining \eqref{eq:lem_eps_1} and  \eqref{eq:lem_eps_2} we get
	$
		L_2(\hat u) \geq \Ec_{\{1:r\}} \geq  L_2(\hat u)
	$
	which completes the proof.
\end{proof}

\subsubsection{Proof of Theorem \ref{thm:update}}\label{app:thm2}
\begin{proof}
The validity of the update scheme \eqref{eq:update_eps_lr} follows directly from 
the derivation described between \eqref{eq:R_z} and \eqref{eq:update_eps_lr}.
It remains to prove the statement on the complexity.
The update \eqref{eq:update_eps_lr} needs only a constant number of Givens rotations
which act on submatrix of size $5 \times 4.$
Thus, the computational cost for an update step \eqref{eq:update_eps_lr} is constant, $O(1),$
and so computing the entire vector  $[\Ec_{\{1:1\}}, \Ec_{\{1:2\}}, \ldots, \Ec_{\{1:N\}}]$ is $O(N).$ 
The memory complexity is $O(N)$ because the matrices $A^{(N)}$ 
is a band-matrix which retains its band-structure after application of the Givens rotations.
\end{proof}

\subsubsection{Proof of Theorem \ref{thm:algorithm}}\label{app:thm3}

\begin{proof}
	We have seen in Section \ref{sec:reduction} that a solution of
		\eqref{eq:pcw_spline_lagrange} can be obtained by solving \eqref{eq:pcw_spline_discrete}  which in turn can be solved by 
	 computing $F_r^*$ for $r=1, \ldots,N$, by \eqref{eq:recurrencePenalized} or its pruned version \eqref{eq:recurrencePenalized_pruning}. 
	It remains to show the assertions on the complexity.
	 	By Theorem \ref{thm:update}, computation of $\Ec_{\{l:r\}}$ for all $l=1,\ldots,r$
	is $O(r)$ time and memory, so computing $F_r^*$ is $O(r)$ time and memory,
	 provided that $F_{l}^*$ has been stored for $l=1,\ldots, r-1.$
	 Hence, computing the final minimum $F_N^*$ is at most $O(N^2)$ time.
	 The final CSSD is computed via determining the minimizing argument of \eqref{eq:epsilon_I} between  two discontinuity locations. This can be done using Reinsch's algorithm 
	 which is in total $O(N)$ time and $O(N)$ memory, see \cite{reinsch1967smoothing} and \cite{de1978practical}.
	 Storing the $F^*$ values needs $O(N)$ memory,
	 and storing the optimal discontinuity locations is $O(N)$ memory as well.
	 The values $\Ec_{\{l:r\}},$ $l=1,\ldots,r,$ can be discarded once $F^{*}_r$ 
	 has been determined.
	 Hence, the overall memory complexity is $O(N).$
\end{proof}

\subsubsection{Proof of Theorem \ref{thm:uniqueness}}\label{app:thm5}

\begin{proof}
	
	(1) For convenience, we first formulate the functionals in \eqref{eq:pcw_spline_discrete} with an explicit dependence on the data $y \in \R^N$:
	\[ 
	\Gc(J, y) = \sum_{I \in \Pc(J)} \Ec_{I,y} + \gamma |J|,
	\]
	where $ \Ec_{I,y}$ is given by \eqref{eq:epsilon_I}, i.e.,	
	\begin{equation*}
		\Ec_{I,y} =  \min_{f \in C^2(I)} p \sum_{i: x_i \in I} \left(\frac{y_i - f(x_i)}{\delta_i}\right)^2   +  (1-p) \int_{I}   (f''(t))^2 \,dt.
	\end{equation*}
	For any interval $I,$ the minimizer of \eqref{eq:epsilon_I} is the corresponding smoothing spline $f_I^\ast$ which is uniquely determined;  cf. e.g. \cite{de1978practical}. Hence, the mapping from data $y_I$ on (a discrete subset of) the interval $I$ to the corresponding smoothing spline on $I$ is well defined. We denote the mapping restricted to data sites in $I$ by 
	\[
	  S_I (y_I) :=  f_I^\ast |_{x_1,\ldots,x_N}.
	\]   
	The smoothing spline solution operator $S_I$ is a linear mapping acting on a linear space of dimension $d_I:= \# \{ x_i: x_i \in I\},$ the number of data sites in $I.$ If and only if  $d_I=2,$ $S_I$ equals the identity $E,$ thus 
	\[
	S_I - E = 0      \quad \Leftrightarrow \quad   d_I=2.
	\]
	Given a (centered midpoint) jump set $J$ with corresponding partitioning $\Pc(J)$ we consider the corresponding block diagonal matrix 
	\begin{equation}\label{eq:Matrix}
		S_{\Pc(J)} = 
		\begin{bmatrix}
				S_{I_1} & 0 & \cdots & 0\\
				0  & 	S_{I_2} &  \cdots & 0 \\			
				\vdots & \ddots & \ddots & \vdots  \\			
				0 & \cdots  & 0  &	 S_{I_k}
			\end{bmatrix}
	\end{equation}
	where each $I \in \Pc(J)$ corresponds to a diagonal block containing the (matrix representation) of the spline solution operator $S_I.$
	Here $k$ denotes the number of parts in $\Pc(J).$
	$S_{\Pc(J)}$ represents the solution operator of \eqref{eq:pcw_spline_discrete} for any given $J.$ 
	We make the straight-forward but important observation that, for (centered midpoint) jump sets $J,J'$  
	\begin{equation}\label{eq:dIcharacterize}
		S_{\Pc(J)} = S_{\Pc(J')}  \quad \Leftrightarrow \quad \text{For all $I \in \Pc(J)$ with $d_I > 2$, $I \in \Pc(J'),$}	
	\end{equation}
	i.e., those segments of the corresponding partitions which have length at least three need to belong to both partitions.
	This leads to the following nomenclature: We call two (midpoint) jump sets $J,J'$ or partitions $\Pc,\Pc'$ {\em essentially different}, if their corresponding matrices $S_{\Pc(J)}$ are not equal.
	
	(2) Next, we aim at showing that for  two essentially different jump sets $J,J'$ we have 
	\begin{equation}\label{eq:ZeroSet}
	\Nc(J,J') = \{ y \in \R^N: \Gc(J, y) = \Gc(J', y) \} 
	\end{equation}
	has Lebesgue measure zero.
	Actually, we show that the Lebesgue measure of the level sets 
	\begin{equation}\label{eq:NiveauSet}
		L_{J,c} =  \{ y \in \R^N: \Gc(J, y) = c \}, \quad c \in \R,
	\end{equation}
    equals $0$ which in turn implies  \eqref{eq:ZeroSet}.
	To this end we plug in the smoothing spline solution operator $S_{\Pc(J)}$ to get 
	\begin{equation}
		\Ec_{I,y} =   p \sum_{i: x_i \in I} \left(\frac{y_i - S_{I}y(x_i)}{\delta_i}\right)^2   +  (1-p) \int_{I}   (S_{I}y''(t))^2 \,dt.
	\end{equation} 
	Please note that by a slight abuse of notation we also use the symbol $S_{I}$ to denote the smoothing spline operator (without sampling at the data sites) here. (The related identification mapping is given by spline interpolation.) 
	If $d_I \leq 2,$ $\Ec_{I,y}$ equals $0$ for all  data $y.$
	It is well-known (cf. \cite{de1978practical}) that the integral on the right-hand side can be written as a quadratic functional in the data $y,$ i.e.,
	\begin{equation}
	   \int_{I}   (S_{I}y''(t))^2 \,dt = y^T  A  y = y^T \Delta^T K \Delta y,
	\end{equation} 
	with a symmetric square matrix $A$ of dimension $d_I,$ 
	a second difference matrix $\Delta$ and an invertible square matrix $K$ of dimension $d_I-2.$
		(Details may be found for instance in \cite{de1978practical}.)	
	The kernel of $\Delta, $ and thus $A$ equals the subspace of affine-linear functions, or more precisely, samples thereof. 
The gradient of the quadratic form $y \mapsto  y^T  A  y$
	equals $2\Delta^T K \Delta y$ which is non-zero whenever $y$ does not equal (the sample of) an affine linear function. 
	Hence, choosing a complement of the two-dimensional subspace of affine-linear functions, the corresponding level sets 
	$\{ y: y^TAy =\mathrm{const} \}$ are Lebesgue zero sets. In turn,  $\{ y\in \R^{d_I}: y^TAy = \mathrm{const} \}$ as well as	
	$\{ y\in \R^{d_I}: \Gc(J, y) = \mathrm{const} \}$ are  Lebesgue zero sets.
	This shows \eqref{eq:NiveauSet}, and in consequence \eqref{eq:ZeroSet}.
	 
	(3) Now we may use the statement of \eqref{eq:ZeroSet} to show the assertions of  Theorem~\ref{thm:uniqueness}.
	We first show that  
	the minimizing values $f(x_i)$ at the data sites $x_i$ are uniquely determined  for a.e. input $y \in \R^N.$
	To see this, we exclude those $y' \in \R^N$ for which there are at least two essentially different  $J,J'$ such that $y \in	\Nc(J,J').$
	Let us call this exclusion set $Y'.$
	Since all $	\Nc(J,J')$ are sets of Lebesgue measure zero and there are only finitely many different midpoint jump sets, the corresponding union w.r.t. $J,J'$ is a set of Lebesgue measure zero as well. Hence, $Y'$ is a zero set, and on its complement the minimizing essential partition is unique. 
 	Next, we show that minimizing segments $I^\ast$ of size at least three are uniquely determined for  a.e. $y \in \R^N.$
 	To this end, we exclude the same zero set $Y'$ as above. On the complement of $Y',$ the essential minimizing jump set is unique. Further, for any two $J,J'$  with  $S_{\Pc(J)} = S_{\Pc(J')}$  (equivalent jump sets), both corresponding partitions need to contain the same intervals of size larger than two by \eqref{eq:dIcharacterize}. Thus, they are uniquely determined on $Y'.$
 	Consequently, by the uniqueness of the smoothing spline,  the minimizing function $f$ is unique in $ [ x_{\min}, x_{\max} ]$ 
	where $x_{\min}, x_{\max}$ 	denote the minimal and maximal argument value of the $x_i$ in the minimizing segment $I^\ast$ of length at least three.
	If all segments are of length at least three, for any two $J,J'$  with  $S_{\Pc(J)} = S_{\Pc(J')}$ implies $J=J'$ (midpoint centered jump sets.) 
	Hence, a corresponding jump set $J$ uniquely represents its matrix in $S_{\Pc(J)}$ in \eqref{eq:Matrix} (actually its equivalence class).
	In consequence, on $Y',$ the segments taking the minimal value in \eqref{eq:pcw_spline_discrete} unique. The corresponding statement for the function $f$ is again a consequence of the uniqueness of the smoothing spline.
	Together, this shows all assertions of the theorem and completes the proof.
\end{proof}

\newpage
\subsection{Compatibility with pruning strategies}
To use the FPVI pruning, we reformulate \eqref{eq:recurrencePenalized} as follows:
\begin{equation}\label{eq:recurrencePenalized_pruning}
    F^*_r 
    = \min \Big\{\Ec_{\{ 1:r\}}; \min_{l= r-1, r-2, \ldots, 2} \Ec_{\{ l:r\}} + \gamma + F^*_{l-1}\Big\}.
\end{equation}
The first term in the parenthesis  captures the cost of using a classical spline on the entire interval.
The second term captures the cost of concatenating a continuous spline on $\{l: r\}$ with the already computed optimal discontinuous spline on 
$\{1:l-1\}.$  
We compute $F^*_r$ as follows. We define $F^*_{r, r} := \Ec_{\{ 1:r\}}$ 
and sequentially compute,  for $l=r-1, r-2, \ldots, 2,$ 
\begin{equation}\label{eq:recurrencePenalized_pruning_2}
    F^*_{l,r} 
    = \min \Big\{F^*_{l+1,r}; \Ec_{\{ l:r\}} + \gamma + F^*_{l-1} \Big\}.
\end{equation}
Then, $F^*_r = F^*_{2,r}.$
We observe that $F^*_{l} \geq 0$ for all $l = 1, \ldots, r$ and
  $\Ec_{\{l:r\}}$ is non-decreasing for the reverse order $l = r-1, r-2 \ldots, 2.$
Thus,  if 
\begin{equation}\label{eq:pruning_condition}
\Ec_{\{l:r\}} + \gamma \geq F^*_{l, r},
\end{equation}
we obtain, for all  $l' \leq l,$ 
\[
F^*_{l, r} \leq \Ec_{\{l:r\}} + \gamma \leq  \Ec_{\{l':r\}} + \gamma \leq \Ec_{\{l':r\}} + \gamma + F^*_{l'} \leq F^*_{l', r}.
\]
So when condition \eqref{eq:pruning_condition} is fulfilled we may break the $l$-loop
and set $F^*_r = F^*_{l,r}.$ 

It remains to show that the pruning is compatible with the order of computation of the $\Ec_{l:r}.$
This can be seen as follows:
The first term in parenthesis $\Ec_{\{ 1:r\}}$  of \eqref{eq:recurrencePenalized_pruning}
is precomputed for $r = 1, \ldots, N,$ in that order.
By Theorem~\ref{thm:update}, this requires linear time  and memory, and thus is  uncritical for the overall complexity.
To use the computation scheme \eqref{eq:recurrencePenalized_pruning_2},
we need to compute the $\Ec_{\{ l:r\}}$ in reverse order,
meaning in the order $l = r-1, r-2, \ldots, 2.$
But this can simply be achieved by applying the update scheme on the flipped data vector.

When using the PELT pruning, the minimum 
in \eqref{eq:recurrencePenalized} has to be  sought over an active set of left bounds $ \Ac_r,$ so 
 $F^*_r 
    =  \min_{l \in \Ac_r} \Big\{  \Ec_{\{ l:r\}} + \gamma + F^*_{l-1} \Big\}.$
    After step $r,$ the active set is pruned according to $\Ac_{r+1} = \{ l \in \Ac_r \cup \{r\}: F^*_{l-1} + \Ec_{\{ l:r\}} \leq F^*_r\}.$
As $\Ac_{r+1} \subset \Ac_r \cup \{r\},$
we may compute the splines energies in the update order
$\Ec_{\{l,r\}} \to \Ec_{\{l,r+1\}}$ for all active left bounds $l \in \Ac_{r+1}$ with $l < r;$
and if $l=r$ we start a new recursion by computing the initial state 
corresponding to $\Ec_{\{r:r+1\}}.$
As for FPVI, the case with no discontinuities can be handled by initializing $F^*_r = \Ec_{\{ 1:r\}}$ for all $r = 1, \ldots, N,$
so that the recursion can be started with $r = 3$ and $\Ac_{3} = \{2\}.$

\subsection{Additional Experiments}

\subsubsection{Experiment on vector-valued data}

\begin{figure}[!t]
\centering
	\includegraphics[width=0.9\textwidth]{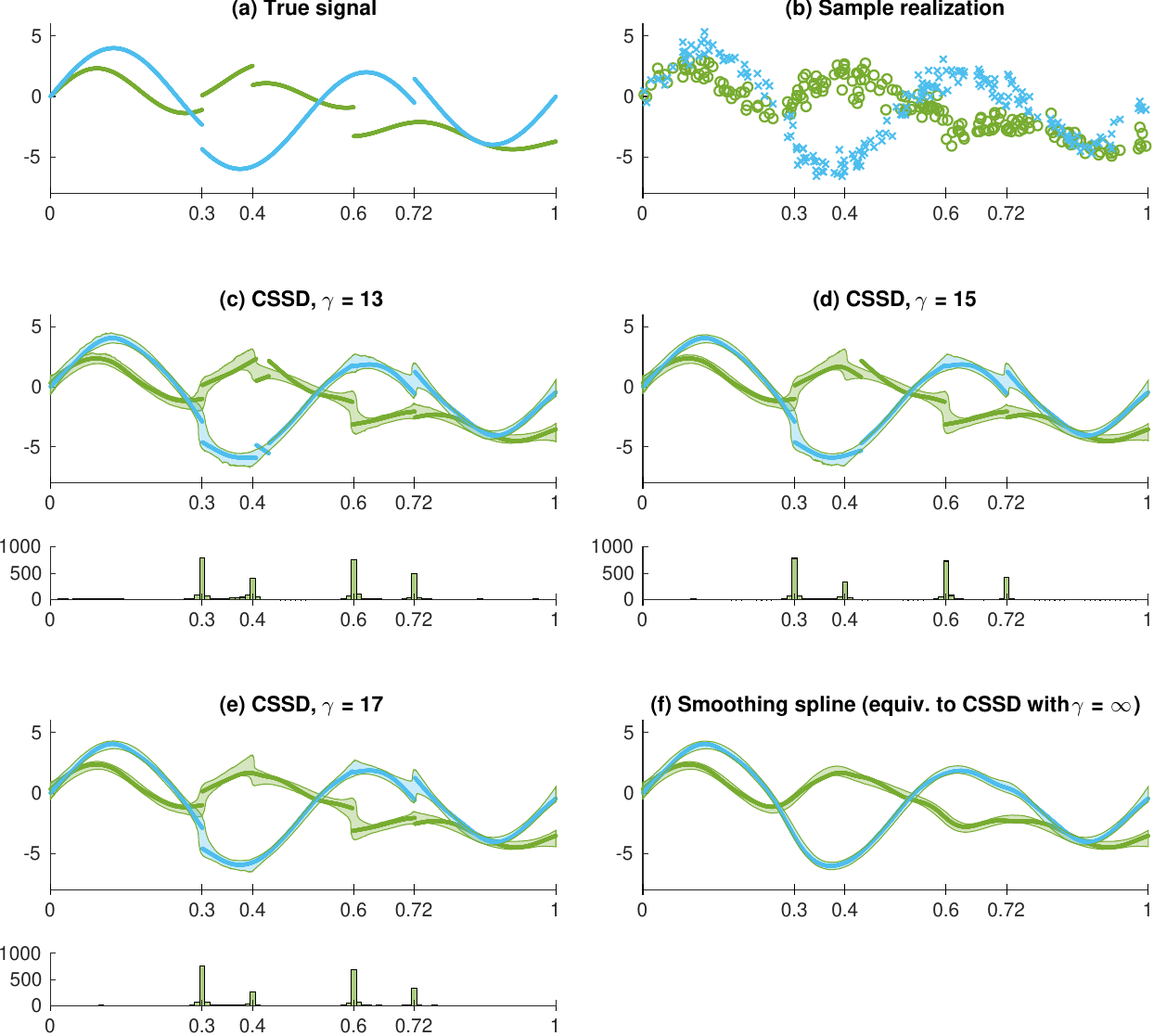}
	\caption{Here a y-value consists of two elements per data site $x_i$, so $y_i \in \R^2$: 
	The first component  are noisy samples of the signal used in Figure~\ref{fig:synthetic},  rescaled by the factor 4, 
	and the second component are noisy samples of the HeaviSine signal where $N = 200,$ $\sigma = 0.6.$
		The description of the subplots matches that of Figure~\ref{fig:synthetic}. Here $p=0.9999$ was used. 
	}
	\label{fig:vector}
\end{figure}
Figure~\ref{fig:vector} shows a vector-valued variant of the last experiment.
The data is generated according to $y_i = [4g_1(x_i), g_2(x_i)] + \epsilon_i \in \R^2,$ for $i=1, \ldots, N,$ where $\epsilon_{i}$ is a zero mean Gaussian random vector with $\sigma = 0.6$ in both components, and $x$ consists of $N = 200$ random values uniformly distributed in the interval $[0,1].$
The highest detection rate for discontinuities is observed at the shared discontinuity $x = 0.3,$ 
and at the discontinuities $x = 0.6$ of $g_1$ and $x =0.72$ of $g_2,$ which have a relatively large jump height. The discontinuity at $x =0.4$ of $g_1$ has the lowest detection rate. This is because jump height is  relatively small and because the signal is smooth in the other component $g_2$.

\subsubsection{Detailed Runtimes for the Geyser data}
The following table shows the runtimes for
computing the CSSD solutions for 
the $(p,\gamma)$ values sampled by the optimizer for
determination of a hyperparameter using cross validation for the  Geyser data.
It illustrates the dependency of the runtime on the hyperparameters.

\begin{center}
\begin{longtable}{rrrrr}
	\toprule 
	\multicolumn{2}{c}{Parameters} & \multicolumn{3}{c}{Runtimes in seconds} \\ 
	$p$ & $\gamma$ & Prop. + PELT & Prop. + FPVI & Baseline \\ 
	\toprule 
	0.590 & 526.700 & 0.367 & 0.173 & 60.084 \\
0.114 & 527.176 & 0.012 & 0.010 & 59.848 \\
0.646 & 0.429 & 0.085 & 0.066 & 1.419 \\
0.493 & 0.658 & 0.077 & 0.069 & 1.510 \\
0.103 & 527.187 & 0.011 & 0.010 & 59.577 \\
0.513 & 4831.609 & 0.011 & 0.010 & 59.585 \\
\bottomrule
\end{longtable}	
\end{center}

\end{document}